\documentclass[submission,copyright,creativecommons,sharealike]{eptcs}
\providecommand{\event}{ACT 2022} %
\providecommand{\thisvolume}[1]{this volume of EPTCS, Open Publishing Association}
\usepackage{breakurl}             %
\usepackage{underscore}           %

\usepackage[utf8]{inputenc}
\usepackage[T1]{fontenc} %

\usepackage{microtype}
\usepackage{amsfonts}
\usepackage{amssymb}
\usepackage{amsthm}
\usepackage{amsmath}
\usepackage{caption}
\usepackage{etoolbox}
\usepackage{stmaryrd}
\usepackage{ifthen}
\usepackage{suffix}
\usepackage{verbatim} %

\usepackage{xfrac}

\usepackage{bm}
\usepackage{physics}
\usepackage[safe]{tipa} %
\usepackage{mathtools} %
\usepackage{color}
\usepackage{wrapfig}
\usepackage{tikz}
\usepackage{tikzit}
\usetikzlibrary{cd, babel}

\usepackage[framemethod=tikz]{mdframed}

\usepackage[backend=biber,style=numeric-comp,sorting=none,date=short]{biblatex}
\appto{\bibsetup}{\sloppy}

\newcommand*{\citep}[1]{\parencite{#1}}
\newcommand*{\citet}[1]{\textcite{#1}}

\usepackage[capitalize]{cleveref} %

\usepackage{graphicx}

\usepackage{pict2e}

\def\mdash{{\hbox{-}}}

\makeatletter
\newcommand{\adjunction}{\@ifstar\named@adjunction\normal@adjunction}
\newcommand{\normal@adjunction}[4]{%
  #1\colon #2%
  \mathrel{\vcenter{%
    \offinterlineskip\m@th
    \ialign{%
      \hfil$##$\hfil\cr
      \longrightharpoonup\cr
      \noalign{\kern-.3ex}
      \smallbot\cr
      \longleftharpoondown\cr
    }%
  }}%
  #3 \noloc #4%
}
\newcommand{\named@adjunction}[4]{%
  #2%
  \mathrel{\vcenter{%
    \offinterlineskip\m@th
    \ialign{%
      \hfil$##$\hfil\cr
      \scriptstyle#1\cr
      \noalign{\kern.1ex}
      \longrightharpoonup\cr
      \noalign{\kern-.3ex}
      \smallbot\cr
      \longleftharpoondown\cr
      \scriptstyle#4\cr
    }%
  }}%
  #3%
}
\newcommand{\longrightharpoonup}{\relbar\joinrel\rightharpoonup}
\newcommand{\longleftharpoondown}{\leftharpoondown\joinrel\relbar}
\newcommand\noloc{%
  \nobreak
  \mspace{6mu plus 1mu}
  {:}
  \nonscript\mkern-\thinmuskip
  \mathpunct{}
  \mspace{2mu}
}
\newcommand{\smallbot}{%
  \begingroup\setlength\unitlength{.15em}%
  \begin{picture}(1,1)
  \roundcap
  \polyline(0,0)(1,0)
  \polyline(0.5,0)(0.5,1)
  \end{picture}%
  \endgroup
}
\makeatother

\let\op=\relax

\let\innerprod=\relax
\newcommand{\innerprod}[2]{\left\langle #1 {,} \: #2 \right\rangle}

\DeclareMathOperator{\op}{^\text{op}}

\newcommand{\nn}{{\mathbb{N}}}
\newcommand{\rr}{{\mathbb{R}}}
\newcommand{\Tt}{{\mathbb{T}}}

\newcommand{\Cat}[1]{\mathbf{#1}}
\newcommand{\cat}[1]{\mathcal{#1}}
\newcommand{\Fun}[1]{\mathsf{#1}}

\newcommand{\Kl}{\mathcal{K}\mspace{-2mu}\ell}

\DeclareMathOperator*{\E}{\mathbb{E}}

\renewcommand{\d}{\mathrm{d}}

\newcommand{\Da}{{\mathcal{D}}}

\newcommand{\Fa}{{\mathcal{F}}}

\newcommand{\Pa}{{\mathcal{P}}}

\DeclareMathOperator{\id}{\mathsf{id}}

\DeclareMathOperator{\Set}{\Cat{Set}}

\newcommand{\xto}[1]{\xrightarrow{#1}}

\makeatletter
\newcommand{\mathoverlap}[2]{\mathpalette\mathoverlap@{{#1}{#2}}}
\newcommand{\mathoverlap@}[2]{\mathoverlap@@{#1}#2}
\newcommand{\mathoverlap@@}[3]{\ooalign{$\m@th#1#2$\crcr\hidewidth$\m@th#1#3$\hidewidth}}
\makeatother

\newcommand{\klcirc}{\bullet} %
\newcommand*{\smallklcirc}{\raisebox{0.18ex}{\scalebox{0.66}{$\klcirc$}}}
\newcommand{\klto}{\mathoverlap{\rightarrow}{\smallklcirc\,}}

\makeatletter
\providecommand*{\xmapstofill@}{%
  \arrowfill@{\mapstochar\relbar}\relbar\rightarrow
}
\providecommand*{\xmapsto}[2][]{%
  \ext@arrow 0395\xmapstofill@{#1}{#2}%
}

\makeatletter
\def\slashedarrowfill@#1#2#3#4#5{%
  $\m@th\thickmuskip0mu\medmuskip\thickmuskip\thinmuskip\thickmuskip
   \relax#5#1\mkern-7mu%
   \cleaders\hbox{$#5\mkern-2mu#2\mkern-2mu$}\hfill
   \mathclap{#3}\mathclap{#2}%
   \cleaders\hbox{$#5\mkern-2mu#2\mkern-2mu$}\hfill
   \mkern-7mu#4$%
}
\def\rightslashedarrowfill@{%
  \slashedarrowfill@\relbar\relbar\mapstochar\rightarrow}
\newcommand\xslashedrightarrow[2][]{%
  \ext@arrow 0055{\rightslashedarrowfill@}{#1}{#2}}
\makeatother

\theoremstyle{definition}
\newtheorem{defn}{Definition}[section]

\newtheorem{ex}[defn]{Example}

\newtheorem{rmk}[defn]{Remark}

\newtheorem*{rmk*}{Remark}

\newtheorem{prop}[defn]{Proposition}
\newtheorem{prop*}{Proposition}

\newtheorem{lemma}[defn]{Lemma}
\newtheorem{thm}[defn]{Theorem}

\newtheorem*{thm*}{Theorem}
\newtheorem*{cor*}{Corollary}

\theoremstyle{remark}
\newtheorem{nproof}[defn]{Proof}

\definecolor{darkblue}{rgb}{0,0,0.7}

\usepackage{enumitem}

\tikzstyle{xshiftu}=[shift = {(#1, 0)}]
\tikzstyle{yshiftu}=[shift = {(0, #1)}]

\tikzstyle{dot}=[inner sep=0.25mm,minimum width=1mm,minimum height=1mm,draw,shape=circle,text depth=-0.2mm]
\tikzstyle{white dot}=[dot,fill=white, draw=black]
\tikzstyle{action}=[dot,fill=white,scale=0.667,inner sep=0.5mm]
\tikzstyle{copier}=[dot,fill=white,scale=2.0]
\tikzstyle{black copier}=[dot,fill=black,scale=2.0]

\tikzstyle{box}=[fill=white, draw=black, shape=rectangle]
\tikzstyle{medium box}=[fill=white, draw=black, shape=rectangle, minimum width=1.5cm, minimum height=0.66cm]
\tikzstyle{arrow box}=[fill=white, draw, shape=rectangle,minimum height=5mm,yshift=-0.5mm,minimum width=5mm]

\tikzstyle{effect}=[regular polygon, regular polygon sides=3,draw]
\tikzstyle{state0}=[regular polygon, regular polygon sides=3,draw,shape border rotate=0]
\tikzstyle{state90}=[regular polygon, regular polygon sides=3,draw,shape border rotate=90]
\tikzstyle{state180}=[regular polygon, regular polygon sides=3,draw,shape border rotate=180]
\tikzstyle{state270}=[regular polygon, regular polygon sides=3,draw,shape border rotate=270]

\tikzstyle{scalar}=[diamond,draw,inner sep=1pt]

\tikzstyle{discarder}=[my ground,draw,inner sep=0pt,minimum width=4.2pt,minimum height=11.2pt,anchor=input,rotate=90]
\tikzstyle{discarder0}=[my ground,draw,inner sep=0pt,minimum width=4.2pt,minimum height=11.2pt,anchor=input,rotate=0]

\tikzstyle{pointy1}=[->]
\tikzstyle{midpoint1}=[-, {postaction={decorate,decoration={markings, mark=at position .5 with {\arrow{>}}}}}]
\tikzstyle{midpointy1pointy}=[->, {postaction={decorate,decoration={markings, mark=at position .5 with {\arrow{>}}}}}]
\tikzstyle{dashed1}=[-, dashed]
\tikzstyle{dotted1}=[-, dotted]
\tikzstyle{dash-pointy}=[->, dashed]

\input{strings.tikzdefs}

\ifexternalizetikz\tikzexternaldisable\fi
\newsavebox\sbground
\savebox\sbground{%
  \begin{tikzpicture}[baseline=0pt]
    \draw (0,-.1ex) to (0,.85ex)
    node[ground IEC,draw,anchor=input,inner sep=0pt,
    minimum width=3.15pt,minimum height=8.4pt,rotate=90] {};
  \end{tikzpicture}%
}
\newcommand{\ground}{\mathord{\usebox\sbground}}

\newsavebox\sbcopier
\savebox\sbcopier{%
  \begin{tikzpicture}[baseline=0pt]
    \node[copier,scale=0.7] (a) at (0,3.8pt) {};
    \draw (a) -- +(-90:.21);
    \draw (a) -- +(45:.21);
    \draw (a) -- +(135:.21);
  \end{tikzpicture}}

\ifexternalizetikz\tikzexternalenable\fi

\newsavebox\bsbcopier
\savebox\bsbcopier{%
  \begin{tikzpicture}[baseline=0pt]
    \node[black copier,scale=0.7] (a) at (0,3.8pt) {};
    \draw (a) -- +(-90:.21);
    \draw (a) -- +(45:.21);
    \draw (a) -- +(135:.21);
  \end{tikzpicture}}
\newcommand{\bcopier}{\mathord{\usebox\bsbcopier}}
\ifexternalizetikz\tikzexternalenable\fi

\usepackage{mathabx}
\usepackage{quiver}

\author{Toby St. Clere Smithe
  \institute{Topos Institute}
  \email{toby@topos.institute}
}
\def\authorrunning{T. St. Clere Smithe}

\hypersetup{
  pdftitle={Dynamical Systems via Generalized Polynomial Coalgebra},
  pdfauthor={Toby St. Clere Smithe}
}

\makeatletter
\newcommand{\ostar}{\mathbin{\mathpalette\make@circled\star}}
\newcommand{\make@circled}[2]{%
  \ooalign{$\m@th#1\smallbigcirc{#1}$\cr\hidewidth$\m@th#1#2$\hidewidth\cr}%
}
\newcommand{\smallbigcirc}[1]{%
  \vcenter{\hbox{\scalebox{0.77778}{$\m@th#1\bigcirc$}}}%
}
\makeatother

\newcommand{\Poly}[1]{\Cat{Poly}_{\cat{#1}}}
\newcommand{\PolyM}[1]{\Cat{Poly}_{#1}}

\newcommand{\RDyn}{\Cat{RDyn}}
\newcommand{\RDynT}[1]{\RDyn^{\mathbb{#1}}}
\newcommand{\BDyn}{\Cat{BunDyn}}
\newcommand{\BDynT}[1]{\BDyn^{\mathbb{#1}}}

\newcommand{\pMCoalgT}[3]{\Cat{Coalg}_{#2}^{#3}({#1})}

\newcommand{\DynT}[1]{\Cat{Coalg}_{\id}^{\mathbb{#1}}}

\newcommand{\PMCoalgT}[1]{\pMCoalgT{p}{M}{#1}}
\newcommand{\MCoalgT}[1]{\pMCoalgT{-}{M}{#1}}
\newcommand{\Coalg}{\Cat{Coalg}}

\makeatletter
\newcommand{\Hier}{\Cat{Hier}}
\newcommand{\HierE@nostar}{{\Cat{Hier}|_{\cat{E}}}}
\newcommand{\HierE@star}[1]{{\Cat{Hier}_{#1}|_{\cat{E}}}}
\newcommand{\HierE}{\@ifstar{\HierE@star}{\HierE@nostar}}
\makeatother

\newcommand{\deloop}[1]{\mathbf{B}#1}
\newcommand{\Sum}{\sum\limits}

\def\ltri{\triangleleft}
\def\rtri{\triangleright}

\date{11 May 2022}

\title{Open Dynamical Systems as Coalgebras for Polynomial Functors, with Application to Predictive Processing}
\def\titlerunning{Dynamical Systems via Generalized Polynomial Coalgebra}

\begin{document}

\maketitle

\begin{abstract}
  We present categories of open dynamical systems with general time evolution as
  categories of coalgebras opindexed by polynomial interfaces, and show how this
  extends the coalgebraic framework to capture common scientific applications
  such as ordinary differential equations, open Markov processes, and random
  dynamical systems. We then extend Spivak's operad $\Cat{Org}$ to this setting,
  and construct associated monoidal categories whose morphisms represent
  hierarchical open systems; when their interfaces are simple, these categories
  supply canonical comonoid structures. We exemplify these constructions using
  the `Laplace doctrine', which provides dynamical semantics for active
  inference, and indicate some connections to Bayesian inversion and coalgebraic
  logic.
\end{abstract}

\section{Background}

\subsection{Closed dynamical systems and Markov processes}

In this brief section, we recall a `behavioural' approach to dynamical systems
originally due (we believe) to Lawvere; for a pedagogical account, see
\parencite{Lawvere2009Conceptual}. These systems are `closed' in the sense that
they do not require environmental interaction for their evolution, but they
nonetheless form the starting point for our categories of more open systems.

\begin{defn}
  Let \((\Tt, +, 0)\) be a monoid, representing time. Let \(X : \cat{E}\) be
  some space, called the \textit{state space}.  Then a \textit{closed dynamical
    system} \(\vartheta\) \textit{with state space} \(X\) \textit{and time}
  \(\Tt\) is an action of \(\Tt\) on \(X\).  When \(\Tt\) is also an object of
  \(\cat{E}\), then this amounts to a morphism \(\vartheta : \Tt \times X \to
  X\) (or equivalently, a time-indexed family of \(X\)-endomorphisms,
  \(\vartheta(t) : X \to X\)), such that \(\vartheta(0) = \id_X\) and
  \(\vartheta(s + t) = \vartheta(s) \circ \vartheta(t)\).
\end{defn}

\begin{prop} \label{prop:transition-map}
  When time is discrete, as in the case \(\Tt = \nn\), any dynamical system
  \(\vartheta\) is entirely determined by its action at \(1 : \Tt\). That is,
  letting the state space be \(X\), we have \(\vartheta(t) = \vartheta(1)^{\circ
    t}\) where \(\vartheta(1)^{\circ t}\) means ``compose \(\vartheta(1) : X \to
  X\) with itself \(t\) times''.
\end{prop}

\begin{ex} \label{ex:closed-vector-field}
  Suppose \(X : U \to TU\) is a vector field on \(U\), with a corresponding
  solution (integral curve) \(\chi_x : \rr \to U\) for all \(x : U\); that is,
  \(\chi'(t) = X(\chi_x(t))\) and \(\chi_x(0) = x\). Then letting the point
  \(x\) vary, we obtain a map \(\chi : \rr \times U \to U\). This \(\chi\) is a
  closed dynamical system with state space \(U\) and time \(\rr\).
\end{ex}

\begin{prop} \label{prop:closed-dyn-cat}
  Closed dynamical systems with state spaces in \(\cat{E}\) and time \(\Tt\) are
  the objects of the functor category \(\Cat{Cat}(\deloop{\Tt}, \cat{E})\),
  where \(\deloop{\Tt}\) denotes the delooping of the monoid \(\Tt\).  Morphisms
  of dynamical systems are therefore natural transformations.
\end{prop}

We will also often be interested in dynamical systems whose evolution has
`side-effects', such as the generation (or `mixing') of uncertainty or
randomness. We will largely model such systems as Kleisli maps or coalgebras of
monads modelling these side-effects. In the case of uncertainty, the monads will
be so-called \textit{probability monads}, which we will often denote by
$\Pa$. Such a monad $\Pa : \cat{E}\to\cat{E}$ can often be thought of as taking
each set or space $X : \cat{E}$ to the set (or space) $\Pa X$ of probability
distributions over $X$, and each morphism to the corresponding `pushforwards'
map; the monad multiplication is given by ``averaging out'' uncertainty, and the
unit takes a point to the `Dirac' distribution over it. With these ideas in
mind, we can extend the concepts above to cover Markov chains and Markov
processes.

\begin{ex}[Closed Markov chains and Markov processes]
  A closed \textit{Markov chain} is given by a map \(X \to \Pa X\), where \(\Pa
  : \cat{E} \to \cat{E}\) is a probability monad on \(\cat{E}\); this is
  equivalently a \(\Pa\)-coalgebra with time \(\nn\), and an object in
  \(\Cat{Cat}\big(\deloop{\nn}, \Kl(\Pa)\big)\). With more general time \(\Tt\),
  one obtains closed \textit{Markov processes}: objects in
  \(\Cat{Cat}\big(\deloop{\Tt}, \Kl(\Pa)\big)\). More explicitly, a closed
  Markov process is a time-indexed family of Markov kernels; that is, a morphism
  \(\vartheta : \Tt \times X \to \Pa X\) such that, for all times \(s,t : \Tt\),
  \(\vartheta_{s+t} = \vartheta_s \klcirc \vartheta_t\) as a morphism in
  \(\Kl(\Pa)\). Note that composition \(\klcirc\) in \(\Kl(\Pa)\) is given by
  the Chapman-Kolmogorov equation, so this means that
  \[
  \vartheta_{s+t}(y|x) = \int_{x':X} \vartheta_s(y|x') \, \vartheta_t(\d x'|x) \, .
  \]
\end{ex}

\subsection{Polynomial functors}

We will use \textit{polynomial functors} to model the interfaces of our open
systems, following \textcite{Spivak2021Polynomial}. We will assume these to be
functors $\cat{E}\to\cat{E}$ for a locally Cartesian closed category $\cat{E}$,
but we will typically assume that $\cat{E}$ is furthermore concrete, and often
that it is in fact $\Cat{Set}$.

\begin{defn}
  Let $\cat{E}$ be a locally Cartesian closed category, and denote by $y^A$ the
  representable copresheaf $y^A := \cat{E}(A, -) : \cat{E} \to \cat{E}$. A
  \emph{polynomial functor} $p$ is a coproduct of representable functors,
  written $p := \sum_{i : p(1)} y^{p_i}$, where $p(1) : \cat{E}$ is the indexing
  object. The category of polynomial functors in $\cat{E}$ is the full
  subcategory $\Poly{E} \hookrightarrow [\cat{E}, \cat{E}]$ of the
  \(\cat{E}\)-copresheaf category spanned by coproducts of representables. A
  morphism of polynomials is therefore a natural transformation.
\end{defn}

\begin{rmk} %
  Every polynomial functor $P : \cat{E} \to \cat{E}$ corresponds to a bundle $p
  : E \to B$ in $\cat{E}$, for which $B = P(1)$ and for each $i : P(1)$, the
  fibre $p_i$ is $P(i)$. We will henceforth elide the distinction between a
  copresheaf $P$ and its corresponding bundle $p$, writing $p(1) := B$ and $p[i]
  := p_i$, where $E = \sum_i p[i]$. A natural transformation $f : p \to q$
  between copresheaves therefore corresponds to a map of bundles. In the case of
  polynomials, by the Yoneda lemma, this map is given by a `forwards' map $f_1 :
  p(1) \to q(1)$ and a family of `backwards' maps $f^\# : q[f_1(\mdash)] \to
  p[\mdash]$ indexed by $p(1)$, as in the left diagram below. Given $f : p \to
  q$ and $g : q \to r$, their composite $g \circ f : p \to r$ is as in the right
  diagram below.
  \begin{equation*}
    \begin{tikzcd}
      E & {f^*F} & F \\
      B & B & C
      \arrow["{f^\#}"', from=1-2, to=1-1]
      \arrow[from=1-2, to=1-3]
      \arrow["q", from=1-3, to=2-3]
      \arrow["p"', from=1-1, to=2-1]
      \arrow[from=2-1, to=2-2, Rightarrow, no head]
      \arrow["{f_1}", from=2-2, to=2-3]
      \arrow[from=1-2, to=2-2]
      \arrow["\lrcorner"{anchor=center, pos=0.125}, draw=none, from=1-2, to=2-3]
    \end{tikzcd}
    \qquad\qquad
    \begin{tikzcd}
      E & {f^*g^*G} & G \\
      B & B & D
      \arrow["{(gf)^\#}"', from=1-2, to=1-1]
      \arrow[from=1-2, to=1-3]
      \arrow["r", from=1-3, to=2-3]
      \arrow["p"', from=1-1, to=2-1]
      \arrow[from=2-2, to=2-1, Rightarrow, no head]
      \arrow["{g_1 \circ f_1}", from=2-2, to=2-3]
      \arrow[from=1-2, to=2-2]
      \arrow["\lrcorner"{anchor=center, pos=0.125}, draw=none, from=1-2, to=2-3]
    \end{tikzcd}
  \end{equation*}
  where $(gf)^\#$ is given by the \(p(1)\)-indexed family of composite maps
  $r[g_1(f_1(\mdash))] \xto{f^\ast g^\#} q[f_1(\mdash)] \xto{f^\#} p[\mdash]$.
\end{rmk}

We can interpret the type $p(1)$ to be a set or space of `configurations' or
`positions' of a $p$-shaped system, and each $p[i]$ to be the available `inputs'
or `directions' available to the system when it is in configuration/position
$i$.

We now recall a handful of useful facts about polynomials and their morphisms,
each of which is explained in \textcite{Spivak2021Polynomial} and summarized in
\textcite{Spivak2022reference}.

\begin{prop}
  Polynomial morphisms $p \to y$ correspond to sections $p(1) \to \sum_i p[i]$
  of the corresponding bundle $p$.
\end{prop}

\begin{prop}
  There is an embedding of \(\cat{E}\) into \(\Poly{E}\) given by taking objects
  \(X : \cat{E}\) to the linear polynomials \(Xy : \Poly{E}\) and morphisms \(f
  : X\to Y\) to morphisms \((f, \id_X) : Xy \to Yy\).
\end{prop}

\begin{prop}
  There is a symmetric monoidal structure \((\otimes, y)\) on \(\Poly{E}\) that
  we call \textnormal{tensor}, and which is given on objects by \(p\otimes q :=
  \sum_{i:p(1)}\sum_{j:q(1)} y^{p[i]\times q[j]}\) and on morphisms \(f := (f_1,
  f^\#) : p\to p'\) and \(g := (g_1, g^\#) : q\to q'\) by \(f\otimes g := (f_1
  \times g_1, f^\#\times g^\#)\).
\end{prop}

\begin{prop}
  $(\Poly{E}, \otimes, y)$ is symmetric monoidal closed, with internal hom
  denoted $[{-},{=}]$. Explicitly, we have $[p,q] = \sum_{f:p\to q}
  y^{\sum_{i:p(1)} q[f_1(i)]}$. Given an object $A : \cat{E}$, we have $[Ay, y]
  \cong y^A$.
\end{prop}

\begin{prop}
  The composition of polynomial functors $q \circ p :
  \cat{E}\to\cat{E}\to\cat{E}$ induces a monoidal structure on $\Poly{E}$, which
  we denote $\ltri$, and call `composition' or `substitution'. Its unit is again
  $y$. Famously, $\ltri$-comonoids correspond to categories and their comonoid
  homomorphisms are cofunctors \parencite{Ahman2016Directed}. If $\Tt$ is a
  monoid, then the comonoid structure on $y^\Tt$ corresponds witnesses it as the
  category $\deloop{\Tt}$. Monomials of the form $Sy^S$ can be equipped with a
  canonical comonoid structure witnessing the codiscrete groupoid on $S$.
\end{prop}

\section{Open dynamical systems as polynomial coalgebras}

\subsection{Deterministic systems}

\begin{defn} \label{def:poly-dyn}
  A deterministic open dynamical system with interface $p$, state space $S$ and time $\Tt$ is a polynomial morphism $\beta : Sy^S \to [\Tt y, p]$ such that, for any section $\sigma : p \to y$, the induced morphism
  \begin{gather*}
    Sy^S \xto{\beta} [\Tt y, p] \xto{[\Tt y, \sigma]} [\Tt y, y] \xto{\sim} y^\Tt
  \end{gather*}
  is a comonoid homomorphism.
\end{defn}

To see how such a morphism $\beta$ is like an `open' version of the closed
dynamical systems introduced above, note that by the tensor-hom adjunction,
$\beta$ can equivalently be written with the type $\Tt y \otimes Sy^S \to p$. In
turn, such a morphism corresponds to a pair $(\beta^o, \beta^u)$, where
$\beta^o$ is the component `on positions' with the type $\Tt \times S \to p(1)$,
and $\beta^u$ is the component `on directions' with the type $\sum_{t:\Tt}
\sum_{s:S} p[\beta^o(t,s)] \to S$. We will call the map $\beta^o$ the
\textit{output map}, as it chooses an output position for each state and moment
in time; and we will call the map $\beta^u$ the \textit{update map}, as it takes
a state $s:S$, a quantity of time $t:\Tt$, and an `input' in $p[\beta^o(t,s)]$,
and returns a new state. We might imagine the new state as being given by
evolving the system from $s$ for time $t$, and the input as being given at the
position corresponding to $(s,t)$.

But it is not sufficient to consider merely such pairs $\beta =
(\beta^o,\beta^u)$ to be our open dynamical systems, for we need them to be like
`open' monoid actions: evolving for time $t$ then for time $s$ must be
equivalent to evolving for time $t+s$, given the same inputs. It is fairly easy
to prove the following proposition, whose proof we defer until after
establishing the categories $\Coalg^\Tt(p)$.

\begin{prop} \label{prop:closed-sys-y-comon}
  Comonoid homomorphisms $Sy^S \to y^\Tt$ correspond bijectively with closed
  dynamical systems with state space $S : \cat{E}$, in the sense given by
  functors $\deloop{\Tt}\to\cat{E}$.
\end{prop}

This establishes that seeking such a comonoid homomorphism will give us the
monoid action property that we seek, and so it remains to show that a composite
comonoid homomorphism of the form $[\Tt y, \sigma]\circ \beta$ is a closed
dynamical system with the ``right inputs''. Unwinding this composite, we find
that the condition that it be a comonoid homomorphism corresponds to the
requirement that, for any $t:\Tt$, the \textit{closure} $\beta^\sigma :
\Tt\times S\to S$ of $\beta$ by $\sigma$ given by
\[
\beta^\sigma(t) := S \xto{\beta^o(t)^\ast \sigma} \Sum_{s:S} p[\beta^o(t, s)] \xto{\beta^u} S
\]
constitutes a closed dynamical system on $S$. The idea here is that $\sigma$
gives the `context' in which we can make an open system closed, thereby
formalizing the ``given the same inputs'' requirement above.

With this conceptual framework in mind, we are in a position to render open
dynamical systems on $p$ with time $\Tt$ into a category, which we will denote
by $\Coalg^\Tt(p)$. Its objects will be pairs $(S, \beta)$ with $S:\cat{E}$ and
$\beta$ an open dynamical on $p$ with state space $S$; we will often write these
pairs equivalently as triples $(S, \beta^o, \beta^u)$, making explicit the
output and update maps. Morphisms will be maps of state spaces that commute with
the dynamics:

\begin{prop} \label{prop:poly-dyn}
  Open dynamical systems over \(p\) with time \(\Tt\) form a category, denoted
  \(\Coalg^\Tt(p)\).  Its morphisms are defined as follows.  Let \(\vartheta :=
  (X, \vartheta^o, \vartheta^u)\) and \(\psi := (Y, \psi^o, \psi^u)\) be two
  dynamical systems over \(p\). A morphism \(f : \vartheta \to \psi\) consists
  in a morphism \(f : X \to Y\) such that, for any time \(t : \Tt\) and global
  section \(\sigma : p(1) \to \Sum_{i:p(1)} p[i]\) of \(p\), the following
  naturality squares commute:
  \[\begin{tikzcd}
	X & {\Sum_{x:X} p[\vartheta^o(t, x)]} & X \\
	Y & {\Sum_{y:Y} p[\psi^o(t, y)]} & Y
	\arrow["{\vartheta^o(t)^\ast \sigma}", from=1-1, to=1-2]
	\arrow["{\vartheta^u(t)}", from=1-2, to=1-3]
	\arrow["f"', from=1-1, to=2-1]
	\arrow["f", from=1-3, to=2-3]
	\arrow["{\psi^o(t)^\ast \sigma}"', from=2-1, to=2-2]
	\arrow["{\psi^u(t)}"', from=2-2, to=2-3]
  \end{tikzcd}\]
  The identity morphism \(\id_\vartheta\) on the dynamical system \(\vartheta\)
  is given by the identity morphism \(\id_X\) on its state space
  \(X\). Composition of morphisms of dynamical systems is given by composition
  of the morphisms of the state spaces.
  \begin{proof}
    We need to check unitality and associativity of composition. This amounts to
    checking that the composite naturality squares commute. But this follows
    immediately, since the composite of two commutative diagrams along a common
    edge is again a commutative diagram.
  \end{proof}
\end{prop}

We can alternatively state Proposition \ref{prop:closed-sys-y-comon} as follows,
noting that the polynomial $y$ corresponds to a trivial interface, exposing no
configuration to any environment nor receiving any signals from it:

\begin{prop} \label{prop:closed-sys-in-dyn-y}
  \(\DynT{T}(y)\) is equivalent to the classical category
  \(\Cat{Cat}(\deloop{\Tt}, \cat{E})\) of closed dynamical systems in
  \(\cat{E}\) with time \(\Tt\).
  \begin{proof}
    The trivial interface \(y\) corresponds to the trivial bundle \(\id_1 : 1
    \to 1\). Therefore, a dynamical system over \(y\) consists of a choice of
    state space \(S\) along with a trivial output map \(\vartheta^o = \ground :
    \Tt \times S \to 1\) and a time-indexed update map \(\vartheta^u : \Tt
    \times S \to S\). This therefore has the form of a classical closed
    dynamical system, so it remains to check the monoid action. There is only
    one section of \(\id_1\), which is again \(\id_1\). Pulling this back along
    the unique map \(\vartheta^o(t) : S \to 1\) gives \(\vartheta^o(t)^\ast
    \id_1 = \id_S\). Therefore the requirement that, given any section
    \(\sigma\) of \(y\), the maps \(\vartheta^u \circ \vartheta^o(t)^\ast
    \sigma\) form an action means in turn that so does \(\vartheta^u : \Tt
    \times S \to S\). Since the pullback of the unique section \(\id_1\) along
    the trivial output map \(\vartheta^o(t) = \ground : S \to 1\) of any
    dynamical system in \(\DynT{T}(y)\) is the identity of the corresponding
    state space \(\id_S\), a morphism \(f : (\vartheta(\ast), \vartheta^u,
    \ground) \to (\psi(\ast), \psi^u, \ground)\) in \(\DynT{T}(y)\) amounts
    precisely to a map \(f : \vartheta(\ast) \to \psi(\ast)\) on the state
    spaces in \(\cat{E}\) such that the naturality condition \(f \circ
    \vartheta^u(t) = \psi^u(t) \circ f\) of Proposition
    \ref{prop:closed-dyn-cat} is satisfied, and every morphism in
    \(\Cat{Cat}(\deloop{\Tt}, \cat{E})\) corresponds to a morphism in
    \(\DynT{T}(y)\) in this way.
  \end{proof}
\end{prop}

Now that we know that our concept of open dynamical system subsumes closed
systems, let us consider some more examples.

\begin{ex}
  Consider a dynamical system \((S, \vartheta^o, \vartheta^u)\) with outputs but
  no inputs.  Such a system has a `linear' interface $p := Iy$ for some $I :
  \cat{E}$; alternatively, we can write its interface $p$ as the `bundle' $\id_I
  : I\to I$.  A section of this bundle must again be $\id_I$, and so
  \(\vartheta^o(t)^\ast \id_I = \id_S\). Once again, the update maps collect
  into to a closed dynamical system in \(\Cat{Cat}(\deloop{\Tt}, \cat{E})\);
  just now we have outputs \(\vartheta^o : \Tt \times S \to p(1) = I\) exposed
  to the environment.
\end{ex}

\begin{prop} \label{prop:open-transition-map}
  When time is discrete, as with \(\Tt = \nn\), any open dynamical system \((X,
  \vartheta^o, \vartheta^u)\) over \(p\) is entirely determined by its
  components at \(1 : \Tt\). That is, we have \(\vartheta^o(t) = \vartheta^o(1)
  : X \to p(1)\) and \(\vartheta^u(t) = \vartheta^u(1) : \sum_{x:X}
  p[\vartheta^o(x)] \to X\). A discrete-time open dynamical system is therefore
  a triple \((X, \vartheta^o, \vartheta^u)\), where the two maps have types
  \(\vartheta^o : X \to p(1)\) and \(\vartheta^u : \sum_{x:X} p[\vartheta^o(x)]
  \to X\).
  \begin{proof}
    Suppose \(\sigma\) is a section of \(p\). We require each closure
    \(\vartheta^\sigma\) to satisfy the flow conditions, that
    \(\vartheta^\sigma(0) = \id_X\) and \(\vartheta^\sigma(t+s) =
    \vartheta^\sigma(t) \circ \vartheta^\sigma(s)\). In particular, we must have
    \(\vartheta^\sigma(t+1) = \vartheta^\sigma(t) \circ
    \vartheta^\sigma(1)\). By induction, this means that we must have
    \(\vartheta^\sigma(t) = \vartheta^\sigma(1)^{\circ t}\) (compare Proposition
    \ref{prop:transition-map}). Therefore we must in general have
    \(\vartheta^o(t) = \vartheta^o(1)\) and \(\vartheta^u(t) = \vartheta^u(1)\).
  \end{proof}
\end{prop}

\begin{ex}
  Suppose \(\dot{x} = f(x, a)\) and \(b = g(x)\), with \(f : X \times A \to TX\)
  and \(g : X \to B\). Then, as for the `closed' vector fields of Example
  \ref{ex:closed-vector-field}, this induces an open dynamical system \((X, \int
  f, g) : \Coalg^\rr(By^A)\), where \(\int f : \rr \times X \times A \to X\)
  returns the \((X,A)\)-indexed solutions of \(f\).
\end{ex}

\begin{ex}
  The preceding example is easily extended to the case of a general polynomial
  interface. Suppose similarly that \(\dot{x} = f(x, a_x)\) and \(b = g(x)\),
  now with \(f : \sum_{x:X} p[g(x)] \to TX\) and \(g : X \to p(1)\). Then we
  obtain an open dynamical system \((X, \int f, g) : \DynT{R}(p)\), where now
  \(\int f : \rr \times \sum_{x:X} p[g(x)] \to X\) is the `update' and \(g : X
  \to p(1)\) the `output' map.
\end{ex}

It is quite straightforward to extend the construction of $\Coalg^\Tt(p)$ to an
opindexed category $\Coalg^\Tt$; we unravel this opindexing explicitly in the
appendix (Proposition \ref{prop:poly-dyn-idx}).

\begin{prop}
  $\Coalg^\Tt$ extends to an opindexed category $\Coalg^\Tt : \Poly{E} \to
  \Cat{Cat}$. On objects (polynomials), it returns the categories above. On
  morphisms of polynomials, we simply post-compose: given $\varphi : p\to q$ and
  $\beta : Sy^S \to [\Tt y, p]$, obtain $Sy^S \to [\Tt y, p] \to [\Tt y, q]$ in
  the obvious way.
\end{prop}

At this point, the reader may be wondering in what sense these open dynamical
systems are coalgebras. To see this, observe that a polynomial morphism $Sy^S
\to q$ is equivalently a map $S\to q(S)$: that is to say, a $q$-coalgebra. By
setting $q = [\Tt y, p]$, we see the connection immediately; to make it clear,
in Proposition \ref{prop:dyn-n-coalg}, we spell it out for the case $\Tt = \nn$.

\subsection{Open Markov processes via stochastic polynomials}

Just as coalgebras $S\to pS$ correspond to discrete-time deterministic open
dynamical systems, coalgebras $S\to p\Pa S$ correspond to discrete-time
\textit{stochastic} dynamical systems when $\Pa$ is a probability monad as
introduced above. We have already seen that `closed' Markov chains correspond to
maps $S\to \Pa S$, and that Markov processes in general time correspond to
functors $\deloop{\Tt}\to\Kl(\Pa)$. Our task in this section is therefore to
connect these two perspectives, extending the categories of deterministic
coalgebras $\Coalg^\Tt(p)$.

Working concretely, it is not hard to spot the relevant adjustment. We therefore
make the following definition.

\begin{defn} \label{def:pM-coalg}
  Let \(M : \cat{E} \to \cat{E}\) be a monad on the category \(\cat{E}\), and
  let \(p : \Poly{E}\) be a polynomial in \(\cat{E}\). Let \((\Tt, +, 0)\) be a
  monoid in \(\cat{E}\), representing time. Then a \(pM\)\textit{-coalgebra with
    time} \(\Tt\) consists in a triple \(\vartheta := (S, \vartheta^o,
  \vartheta^u)\) of a \textnormal{state space} \(S : \cat{E}\) and two morphisms
  \(\vartheta^o : \Tt \times S \to p(1)\) and \(\vartheta^u : \sum_{t:\Tt}
  \sum_{s:\mathbb{S}} p[\vartheta^o(t,s)] \to MS\), such that, for any section
  \(\sigma : p(1) \to \sum_{i:p(1)} p[i]\) of \(p\), the maps \(\vartheta^\sigma
  : \Tt \times S \to MS\) given by
  \[
  \Sum_{t:\Tt} S \xto{\vartheta^o(-)^\ast \sigma} \Sum_{t:\Tt} \Sum_{s:S} p[\vartheta^o(-, s)] \xto{\vartheta^u} MS
  \]
  constitute an object in the functor category \(\Cat{Cat}\big(\deloop{\Tt},
  \Kl(T)\big)\), where \(\deloop{\Tt}\) is the delooping of \(\Tt\) and
  \(\Kl(T)\) is the Kleisli category of \(T\). Once more, we call the closed
  system \(\vartheta^\sigma\), induced by a section \(\sigma\) of \(p\), the
  \textnormal{closure} of \(\vartheta\) by \(\sigma\).
\end{defn}

As before, such $pM$-coalgebras form a category; and these categories in turn
are opindexed by polynomials.

\begin{prop} \label{prop:pM-coalg}
  \(pM\)-coalgebras with time \(\Tt\) form a category, denoted
  \(\PMCoalgT{\Tt}\). Its morphisms are defined as follows.  Let \(\vartheta :=
  (X, \vartheta^o, \vartheta^u)\) and \(\psi := (Y, \psi^o, \psi^u)\) be two
  \(pM\)-coalgebras. A morphism \(f : \vartheta \to \psi\) consists in a
  morphism \(f : X \to Y\) such that, for any time \(t : \Tt\) and global
  section \(\sigma : p(1) \to \Sum_{i:p(1)} p[i]\) of \(p\), the following
  naturality squares commute:
  \[\begin{tikzcd}
	X & {\Sum_{x:X} p[\vartheta^o(t, x)]} & MX \\
	Y & {\Sum_{y:Y} p[\psi^o(t, y)]} & MY
	\arrow["{\vartheta^o(t)^\ast \sigma}", from=1-1, to=1-2]
	\arrow["{\vartheta^u(t)}", from=1-2, to=1-3]
	\arrow["f"', from=1-1, to=2-1]
	\arrow["Mf", from=1-3, to=2-3]
	\arrow["{\psi^o(t)^\ast \sigma}"', from=2-1, to=2-2]
	\arrow["{\psi^u(t)}"', from=2-2, to=2-3]
  \end{tikzcd}\]
  The identity morphism \(\id_\vartheta\) on the \(pM\)-coalgebra \(\vartheta\)
  is given by the identity morphism \(\id_X\) on its state space
  \(X\). Composition of morphisms of \(pM\)-coalgebras is given by composition
  of the morphisms of the state spaces.
\end{prop}

\begin{prop} \label{prop:pM-coalg-idx}
  \(\PMCoalgT{\Tt}\) extends to an opindexed category, \(\MCoalgT{\Tt} :
  \Poly{E} \to \Cat{Cat}\). Suppose \(\varphi : p \to q\) is a morphism of
  polynomials.  We define a corresponding functor \(\pMCoalgT{\varphi}{M}{\Tt} :
  \pMCoalgT{p}{M}{\Tt} \to \pMCoalgT{q}{M}{\Tt}\) as follows.  Suppose \((X,
  \vartheta^o, \vartheta^u) : \pMCoalgT{p}{M}{\Tt}\) is an object
  (\(pM\)-coalgebra) in \(\pMCoalgT{p}{M}{\Tt}\). Then
  \(\pMCoalgT{\varphi}{M}{\Tt}(X, \vartheta^o, \vartheta^u)\) is defined as the
  triple \((X, \varphi_1 \circ \vartheta^o, \vartheta^u \circ {\vartheta^o}^\ast
  \varphi^\#) : \pMCoalgT{q}{M}{\Tt}\), where the two maps are explicitly the
  following composites:
  \begin{gather*}
    \Tt \times X \xto{\vartheta^o} p(1) \xto{\varphi_1} q(1) \, ,
    \qquad
    \Sum_{t:\Tt} \Sum_{x:X} q[\varphi_1 \circ \vartheta^o(t, x)] \xto{{\vartheta^o}^\ast \varphi^\#} \Sum_{t:\Tt} \Sum_{x:X} p[\vartheta^o(t, x)] \xto{\vartheta^u} MX \, .
  \end{gather*}
  On morphisms, \(\pMCoalgT{\varphi}{M}{\Tt}(f) : \pMCoalgT{\varphi}{M}{\Tt}(X,
  \vartheta^o, \vartheta^u) \to \pMCoalgT{\varphi}{M}{\Tt}(Y, \psi^o, \psi^u)\)
  is given by the same underlying map \(f : X \to Y\) of state spaces.
\end{prop}

The opindexed category $\PMCoalgT{\Tt}$ clearly generalizes $\Coalg^\Tt$, since
we can always take $M = \id_{\cat{E}}$. Yet these concrete definitions obscure
the more elegant representation of the objects of $\Coalg^\Tt$ as morphisms
$Sy^S \to [\Tt y, p]$. Our task is therefore to find a setting in which a
similar representation is possible; to do so, we generalize $\Poly{E}$ so that
the backwards components of its morphisms may incorporate `side-effects'
modelled by $M$. We will call the corresponding category $\PolyM{M}$, and will
find that instantiating $\Coalg^\Tt$ in $\PolyM{M}$ recovers $\PMCoalgT{\Tt}$.

We begin by recalling that $\Poly{E}$ is equivalent to the category of
Grothendieck lenses for the self-indexing
\parencite{Spivak2020Poly,Spivak2021Polynomial}: $\Poly{E} \cong \int
\cat{E}/-\op$, where the opposite is taken pointwise on each $\cat{E}/B$. We
will define $\PolyM{M}$ by analogy, using the following indexed
category. Suppose $M$ is a commutative monad on $\cat{E}$ and let $\iota$ denote
the identity-on-objects inclusion $\cat{E} \hookrightarrow \Kl(M)$ given on
morphisms by post-composing with the unit $\eta$ of the monad structure. For
ease of exposition in this short paper, we will assume here that $\cat{E} =
\Set$.

\begin{defn}
  Define the indexed category $\cat{E}_M/- : \cat{E}\op \to \Cat{Cat}$ as
  follows. On objects $B : \cat{E}$, we define $\cat{E}_M/B$ to be the full
  subcategory of $\Kl(M)/B$ on those objects $\iota p : E\klto B$ which
  correspond to maps $E \xto{p} B \xto{\eta_B} MB$ in the image of $\iota$. Now
  suppose $f : C\to B$ is a map in $\cat{E}$. We define $\cat{E}_M/f :
  \cat{E}_M/B \to \cat{E}_M/C$ as follows. The functor $\cat{E}_M/f$ takes
  objects $\iota p : E\klto B$ to $\iota(f^*p) : f^*E\klto C$ where $f^*p$ is
  the pullback of $p$ along $f$ in $\cat{E}$, included into $\Kl(M)$ by
  $\iota$.

  To define the action of $\cat{E}_M/f$ on morphisms $\alpha : (E, \iota
  p:E\klto B) \to (F, \iota q:F\klto B)$, note that since we must have $\iota
  q\klcirc \alpha = \iota p$, $\alpha$ must correspond to a family of maps
  $\alpha_x : p[x] \to Mq[x]$ for $x : B$. Then we can define
  $(\cat{E}_M/f)(\alpha)$ pointwise as $(\cat{E}_M/f)(\alpha)_y := \alpha_{f(y)}
  : p[f(y)] \to Mq[f(y)]$ for $y : C$.
\end{defn}

\begin{defn}
  We define $\PolyM{M}$ to be the category of Grothendieck lenses for
  $\cat{E}_M/-$. That is, $\PolyM{M} := \int \cat{E}_M/-\op$, where the opposite
  is again taken pointwise.
\end{defn}

Unwinding this definition, we find that the objects of $\PolyM{M}$ are the same
polynomial functors as constitute the objects of $\Poly{E}$. The morphisms $f :
p \to q$ are pairs $(f_1, f^\#)$, where $f_1 : B\to C$ is a map in $\cat{E}$ and
$f^\#$ is a family of morphisms $q[f_1(x)]\klto p[x]$ in $\Kl(M)$, making the
following diagram commute:

\begin{equation*}
  \begin{tikzcd}
    \sum_{x:B} Mp[x] & \sum_{b:B} q[f_1(x)] & \sum_{y:C} q[y] \\
    B & B & C
    \arrow["{f^\#}"', from=1-2, to=1-1]
    \arrow[from=1-2, to=1-3]
    \arrow["q", from=1-3, to=2-3]
    \arrow["{\eta_B}^* p"', from=1-1, to=2-1]
    \arrow[from=2-1, to=2-2, Rightarrow, no head]
    \arrow["{f_1}", from=2-2, to=2-3]
    \arrow[from=1-2, to=2-2]
    \arrow["\lrcorner"{anchor=center, pos=0.125}, draw=none, from=1-2, to=2-3]
  \end{tikzcd}
\end{equation*}

\begin{rmk}
  Note that the tensor $\otimes$ extends to $\PolyM{M}$: on objects, it is
  defined identically to the tensor on $\Poly{E}$. On morphisms $f := (f_1,
  f^\#) : p\to q$ and $g := (g_1, g^\#) : p' \to q'$, we define the tensor
  $f\otimes g$ to have forwards component $f_1\times g_1$ as before, and the
  backwards components are defined by $(f\otimes g)^\#_{(x, x')} :=
  q[f_1(x)]\times q'[g_1(x')] \to Mp[x]\times Mp'[x'] \to M\bigl(p[x]\times
  p'[x']\bigr)$, where the second arrow is given by the commutativity of the
  monad $M$. On the other hand, we only get an internal hom satisfying the
  adjunction $\PolyM{M}(p\otimes q, r) \cong \PolyM{M}(p, [q,r])$ when the
  backwards components of morphisms $p\otimes q \to r$ are `uncorrelated'
  between $p$ and $q$.
\end{rmk}

\begin{rmk}
  For $\PolyM{M}$ to behave faithfully like the category $\Poly{E}$ of
  polynomial functors and their morphisms, we should want the substitution
  functors $\cat{E}_M/f : \cat{E}_M/C \to \cat{E}_M/B$ to have left and right
  adjoints. Although we do not spell it out here, it is quite straightforward to
  exhibit the left adjoints. On the other hand, writing $f^*$ as shorthand for
  $\cat{E}_M/f$, we can see that a right adjoint only obtains in restricted
  circumstances. Denote the putative right adjoint by $\Pi_f : \cat{E}_M/B \to
  \cat{E}_M/C$, and for $\iota p : E\klto B$ suppose that $(\Pi_f E)[y]$ is
  given by the set of `partial sections' $\sigma : f^{-1}\{y\} \to TE$ of $p$
  over $f^{-1}\{y\}$ as in the commutative diagram:
  \[\begin{tikzcd}
	& {f^{-1}\{y\}} & {\{y\}} \\
	TE & B & C
	\arrow[from=1-2, to=1-3]
	\arrow[from=1-3, to=2-3]
	\arrow[from=1-2, to=2-2]
	\arrow["f", from=2-2, to=2-3]
	\arrow["\lrcorner"{anchor=center, pos=0.125}, draw=none, from=1-2, to=2-3]
	\arrow["{{\eta_B}^*p}", from=2-1, to=2-2]
	\arrow["\sigma"', curve={height=12pt}, from=1-2, to=2-1]
  \end{tikzcd}\]
  Then we would need to exhibit a natural isomorphism $\cat{E}_M/B(f^*D, E)
  \cong \cat{E}_M/C(D, \Pi_f E)$. But this will only obtain when the `backwards'
  components $h^\#_y : D[y]\to M(\Pi_f E)[y]$ are in the image of
  $\iota$---otherwise, it is not generally possible to pull $f^{-1}\{y\}$ out of
  $M$.
\end{rmk}

Despite these restrictions, we do have enough structure at hand to instantiate
$\Coalg^\Tt$ in $\PolyM{M}$. The only piece remaining is the composition product
$\ltri$, but for our purposes it suffices to define its action on objects, which
is identical to its action on objects in $\Poly{E}$\footnote{We leave the full
exposition of $\ltri$ in $\PolyM{M}$ to the forthcoming extended version of this
paper.}, and then consider $\ltri$-comonoids in $\PolyM{M}$. The comonoid laws
force the structure maps to be deterministic (\textit{i.e.}, in the image of
$\iota$), and so $\ltri$-comonoids in $\PolyM{M}$ are just $\ltri$-comonoids in
$\Poly{E}$.

Finally, we note that we can define morphisms $\beta : Sy^S \to [\Tt y, p]$:
these again just correspond to morphisms $\Tt y \otimes Sy^S \to p$, and the
condition that the backwards maps be uncorrelated between $\Tt y$ and $p$ is
satisfied because $\Tt y$ has a trivial exponent. Unwinding such a $\beta$
according to the definition of $\PolyM{M}$ indeed gives precisely a pair
$(\beta^o, \beta^u)$ of the requisite types; and a comonoid homomorphism $Sy^S
\to y^\Tt$ in $\PolyM{M}$ is precisely a functor $\deloop{\Tt} \to \Kl(M)$,
thereby establishing equivalence between the objects of $\Coalg^{\Tt}(p)$
established in $\PolyM{M}$ and the objects of $\PMCoalgT{\Tt}$. The equivalence
between the hom-sets is established by a similar unwinding. All told, in this
section, we have sketched the proof of the following theorem:

\begin{thm}
  Constructing $\Coalg^{\Tt}(p)$ in $\PolyM{M}$ yields a category equivalent to
  $\PMCoalgT{\Tt}$.
\end{thm}

\subsection{Random dynamical systems and bundle systems}

In the analysis of stochastic systems, it is often fruitful to consider two
perspectives: on one side, one considers explicitly the evolution of the
distribution of the states of the system, by following (for instance) a Markov
process, or Fokker-Planck equation. On the other side, one considers the system
as if it were a deterministic system, perturbed by noisy inputs, giving rise to
the frameworks of stochastic differential equations and associated
\textit{random dynamical systems}.

Whereas a (closed) Markov process is typically given by the action of a `time'
monoid on an object in a Kleisli category of a probability monad, a (closed)
random dynamical system is given by a \textit{bundle} of closed dynamical
systems, where the base system is equipped with a probability measure which it
preserves: the idea being that a random dynamical system can be thought of as a
`random' choice of dynamical system on the total space at each moment in time,
with the base measure-preserving system being the source of the randomness
\parencite{Arnold1998Random}.

This idea corresponds in non-dynamical settings to the notion of \textit{randomness pushback} \parencite[Def. 11.19]{Fritz2019synthetic}, by which a stochastic map \(f : A\to\Pa B\) can be presented as a deterministic map \(f^\flat : \Omega\times A\to B\) where \((\Omega,\omega)\) is a probability space such that, for any \(a:A\), pushing \(\omega\) forward through \(f^\flat(\mdash,a)\) gives the state \(f(a)\); that is, \(\omega\) induces a random choice of map \(f^\flat(\omega,\mdash) : A\to B\).
Similarly, under `nice' conditions, random dynamical systems and Markov processes do coincide, although they have different suitability in applications.

In this section, we sketch how the generalized-coalgebraic structures developed
above extend also to random dynamical systems, though with most details deferred
to the Appendix.  By observing that we can also `open up' the base system of a
random dynamical system, we obtain furthermore a notion of \textit{open bundle
  system}: a bundle of dynamical systems that is coherently `open' over
polynomials both in the total space and the base space.

\begin{defn}
  Suppose \(\cat{E}\) is a category equipped with a probability monad \(\Pa :
  \cat{E} \to \cat{E}\) and a terminal object \(1 : \cat{E}\).  A
  \textnormal{probability space} in \(\cat{E}\) is an object of the slice \(1 /
  \Kl(\Pa)\) of the Kleisli category of the probability monad under \(1\).
\end{defn}

\begin{rmk}
  In order to consider polynomials in $\cat{E}$, we will later assume again that
  it is locally Cartesian closed. A simple example of a locally Cartesian closed
  category equipped with a probability monad is the category $\Set$ equipped
  with the monad $\Da$ taking each set to the set of finitely-supported
  probability distributions upon it.
\end{rmk}

\begin{prop}
  There is a forgetful functor \(1/\Kl(\Pa) \to \cat{E}\) taking probability
  spaces \((B, \beta)\) to the underlying spaces \(B\) and their morphisms \(f :
  (A, \alpha) \to (B, \beta)\) to the underlying maps \(f : A \to \Pa B\). We
  will write \(B\) to refer to the space in \(\cat{E}\) underlying a probability
  space \((B, \beta)\), in the image of this forgetful functor.
\end{prop}

\begin{defn} \label{def:metric-sys}
  Let \((B, \beta)\) be a probability space in \(\cat{E}\). A closed
  \textnormal{metric} or \textnormal{measure-preserving} dynamical system \((\vartheta,
  \beta)\) on \((B, \beta)\) with time \(\Tt\) is a closed dynamical system
  \(\vartheta\) with state space \(B : \cat{E}\) such that, for all \(t : \Tt\),
  \(\Pa \vartheta(t) \circ \beta = \beta\); that is, each \(\vartheta(t)\) is a
  \((B, \beta)\)-endomorphism in \(1/\Kl(\Pa)\).
\end{defn}

\begin{prop} \label{prop:metric-sys}
  Closed measure-preserving dynamical systems in \(\cat{E}\) with time \(\Tt\)
  form the objects of a category \(\Cat{Cat}(\deloop{\Tt}, \cat{E})_{\Pa}\)
  whose morphisms \(f : (\vartheta, \alpha) \to (\psi, \beta)\) are maps \(f :
  \vartheta(\ast) \to \psi(\ast)\) in \(\cat{E}\) between the state spaces that
  preserve both flow and measure, as in the following commutative diagram, which
  also indicates their composition:
  
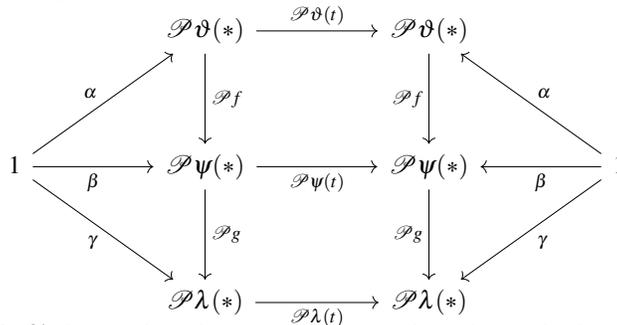
\vspace*{1.75cm}
  \[\begin{tikzcd}[transform canvas={scale=0.875}]
	&& {\Pa\vartheta(\ast)} && {\Pa\vartheta(\ast)} \\
	\\
	1 && {\Pa\psi(\ast)} && {\Pa\psi(\ast)} && 1 \\
	\\
	&& {\Pa\lambda(\ast)} && {\Pa\lambda(\ast)}
	\arrow["\alpha", from=3-1, to=1-3]
	\arrow["\beta"', from=3-1, to=3-3]
	\arrow["\gamma"', from=3-1, to=5-3]
	\arrow["\alpha"', from=3-7, to=1-5]
	\arrow["\beta", from=3-7, to=3-5]
	\arrow["\gamma", from=3-7, to=5-5]
	\arrow["{\Pa\vartheta(t)}", from=1-3, to=1-5]
	\arrow["{\Pa\psi(t)}"', from=3-3, to=3-5]
	\arrow["{\Pa\lambda(t)}"', from=5-3, to=5-5]
	\arrow["{\Pa f}", from=1-3, to=3-3]
	\arrow["{\Pa f}"', from=1-5, to=3-5]
	\arrow["{\Pa g}", from=3-3, to=5-3]
	\arrow["{\Pa g}"', from=3-5, to=5-5]
  \end{tikzcd}\]
  \vspace*{1.25cm}
\end{prop}

\begin{defn}
  Let \((\vartheta, \beta)\) be a closed measure-preserving dynamical system. A
  closed \textnormal{random dynamical system} over \((\vartheta, \beta)\) is an
  object of the slice category \(\Cat{Cat}(\deloop{\Tt}, \cat{E})/\vartheta\);
  it is therefore a bundle of the corresponding functors.
\end{defn}

\begin{ex} \label{ex:brown-sde}
  The solutions \(X(t, \omega; x_0) : \rr_+ \times \Omega \times M \to M\) to a
  stochastic differential equation \(\d X_t = f(t, X_t) \d t + \sigma(t, X_t) \d
  W_t\), where \(W : \rr_+ \times \Omega \to M\) is a Wiener process in \(M\),
  define a random dynamical system \(\rr_+ \times \Omega \times M \to M : (t,
  \omega, x) \mapsto X(t, \omega; x_0)\) over the Wiener base flow \(\theta :
  \rr_+ \times \Omega \to \Omega : (t, \omega) \mapsto W(s+t, \omega) - W(t,
  \omega)\) for any \(s : \rr_+\).
\end{ex}

\begin{defn} \label{def:poly-rdyn}
  Let \((\theta, \beta)\) be a closed measure-preserving dynamical system in
  \(\cat{E}\) with time \(\Tt\), and let \(p : \Poly{E}\) be a polynomial in
  \(\cat{E}\). Write \(\Omega := \theta(\ast)\) for the state space of
  \(\theta\), and let \(\pi : S \to \Omega\) be an object (bundle) in
  \(\cat{E}/\Omega\). An \textnormal{open random dynamical system over} \((\theta,
  \beta)\) \textnormal{on the interface} \(p\) \textnormal{with state space} \(\pi:S \to
  \Omega\) \textnormal{and time} \(\Tt\) consists in a pair of morphisms
  \(\vartheta^o : \Tt \times S \to p(1)\) and \(\vartheta^u : \Sum_{t:\Tt}
  \Sum_{s:S} p[\vartheta^o(t, s)] \to S\), such that, for any global section
  \(\sigma : p(1) \to \Sum_{i:p(1)} p[i]\) of \(p\), the maps \(\vartheta^\sigma
  : \Tt \times S \to S\) defined as
  \[
  \Sum_{t:\Tt} S \xto{\vartheta^o(-)^\ast \sigma} \Sum_{t:\Tt} \Sum_{s:S} p[\vartheta^o(-, s)] \xto{\vartheta^u} S
  \]
  form a closed random dynamical system in \(\Cat{Cat}(\deloop{\Tt},
  \cat{E})/\theta\), in the sense that, for all \(t : \Tt\) and sections
  \(\sigma\), the following diagram commutes:
  
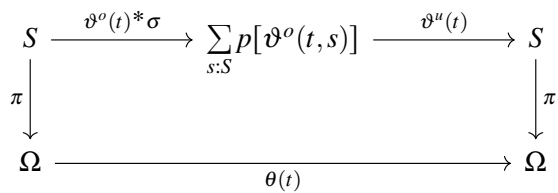
\[\begin{tikzcd}
  S && {\Sum_{s:S} p[\vartheta^o(t, s)]} && S \\
	\Omega &&&& \Omega
	\arrow["\pi"', from=1-1, to=2-1]
	\arrow["\pi", from=1-5, to=2-5]
	\arrow["{\theta(t)}"', from=2-1, to=2-5]
	\arrow["{\vartheta^o(t)^\ast \sigma}", from=1-1, to=1-3]
	\arrow["{\vartheta^u(t)}", from=1-3, to=1-5]
  \end{tikzcd}\]
\end{defn}

\begin{prop} \label{prop:poly-rdyn}
  Let \((\theta, \beta)\) be a closed measure-preserving dynamical system in
  \(\cat{E}\) with time \(\Tt\), and let \(p : \Poly{E}\) be a polynomial in
  \(\cat{E}\). Open random dynamical systems over \((\theta, \beta)\) on the
  interface \(p\) form the objects of a category \(\RDynT{T}(p, \theta)\). See
  Definition \ref{def:cat-open-rds} in the Appendix for details.
\end{prop}

\begin{prop} \label{prop:poly-rdyn-idx}
  The categories \(\RDynT{T}(p, \theta)\) collect into a doubly-indexed category
  of the form \(\RDynT{T} : \Poly{E} \times \Cat{Cat}(\deloop{\Tt},
  \cat{E})_{\Pa} \to \Cat{Cat}\). See Proposition
  \ref{prop:poly-rdyn-idx-details} in the Appendix for details.
\end{prop}

By allowing the base systems of open random dynamical systems instead to be
arbitrary dynamical systems, and then by opening them up similarly, one obtains
notions of \textit{open bundle dynamical system}, and correspondingly
doubly-opindexed categories over pairs of polynomials. Representing these
categories concisely, as we did for the categories $\PMCoalgT{\Tt}$, is the
subject of on-going work, and so we defer the details to the Appendix, in
Definition \ref{defn:bdyn-pbth}, and Propositions \ref{prop:bdyn-pbth-cat},
\ref{prop:bdyn-bth-idx}, and \ref{prop:bdyn-b-idx}.

\section{Hierarchical systems via generalized Org}

In order to exhibit the main example of this paper, we will need to construct,
from the opindexed categories of $p\Pa$-coalgebras introduced above, monoidal
categories whose objects represent the interfaces of hierarchical systems and
whose morphisms represent the hierarchical systems themselves. Informally put,
we will think of a morphism $p\to q$ in such a category as ``a $q$-shaped system
with a $p$-shaped hole''. In order to achieve this, we will in turn adopt and
generalize the operad $\Cat{Org}$ introduced by \textcite{Spivak2021Learners}.

\begin{defn}[Following {\textcite[Def. 2.19]{Spivak2021Learners}}]
  We define a (category-enriched, symmetric, coloured) operad, $\Cat{Org}^\Tt_M$. Its
  objects are polynomials, and for any tuple of polynomials $(p_1, \dots, p_k;
  p')$ of at least length 2, the hom category $\Cat{Org}^\Tt_M(p_1, \dots, p_k;
  p')$ is given by $\pMCoalgT{[p_1 \otimes \dots \otimes p_k, p']}{M}{\Tt}$.
  Note that $\{y \to [\Tt y, [p, p]]\} \cong \{\Tt y \to [p, p]\}$. On any given
  interface $p$, the identity coalgebra is therefore given by the morphism $\Tt
  y \to [p, p]$ that constantly outputs $\id_p$ and has trivial backwards
  component. To define composition, we use the canonical maps
  $[p,q]\otimes[q,r]\to[p,r]$ and $[p,q]\otimes[p',q']\to[p\otimes p', q\otimes
    q']$, the pseudofunctoriality of $\MCoalgT{\Tt}$, and the laxators
  $\pMCoalgT{p}{M}{\Tt}\times \pMCoalgT{q}{M}{\Tt} \to \pMCoalgT{p\otimes
    q}{M}{\Tt}$; since each of these components is associative and unital, the
  composition is well-defined.
\end{defn}

\begin{rmk}
  Spivak's original definition of $\Cat{Org}$ corresponds to the case where $M =
  \id_{\cat{E}}$ and $\Tt = \nn$.
\end{rmk}

For our present purposes, all that is required is to obtain from $\Cat{Org}$ a
(monoidal) (bi)category\footnote{and in fact, we won't even really need to make
use of the monoidal or bicategorical structures here!}. We therefore restrict
$\Cat{Org}^\Tt_M$ to a bicategory $\Hier$ whose objects are again polynomials
and whose hom-categories from $p$ to $q$ are given by $\Cat{Org}^\Tt_M(p, q)$;
it inherits a monoidal structure from the monoidal category associated to the
symmetric operad $\Cat{Org}^\Tt_M$. We will write $\HierE$ to denote the
restriction of $\Hier$ to the linear polynomials $Ay$.

To bring things a little down to earth, first consider a general system $\beta :
p\to q$ in $\Hier$. Recall that $[p, q] = \sum_{f:p\to q} y^{\sum_{i:p(1)}
  q[f_1(i)]}$. $\beta$ is therefore given by a choice of state space $X$ along
with a pair of maps $\beta^o : \Tt\times X \to \PolyM{M}(p,q)$ and $\beta^u :
\sum_{t:\Tt}\sum_{x:X}\sum_{i:p(1)} q[\beta^o(t,s)_1(i)]\to MX$.

To make this a little more comprehensible again, suppose $p = Ay^S$ and $q =
By^T$. Then $\PolyM{M}(p,q) = \cat{E}(A,B)\times \cat{E}(A\times T, S)$, and so
by the universal property of the product, $\beta^o$ is equivalently given by a
pair of maps: a `forwards' output map $\beta^o_1 : \Tt\times X\times A \to B$
and a `backwards' output map $\beta^o_2 : \Tt\times X \times A\times T\to S$; if
this reminds you of a category of lenses, then this is no surprise: the
subcategory of $\Poly{E}$ on the monomials $Ay^S$ is indeed the category of
bimorphic lenses in $\cat{E}$. Finally, the update map simplifies to $\beta^u :
\Tt\times X\times A\times T \to MX$, which updates the state given `forwards'
inputs in $A$ and `backwards' inputs in $T$. We might denote the subcategory of
$\Hier$ on such linear polynomials as $\Cat{HiBi}$, to indicate `hierarchical
bidirectional' systems.

Taking one further step down the ladder of complexity, we briefly consider
systems $\beta : Ay\to By$ in $\HierE$: these are just hierarchical
bidirectional systems where $S = T = 1$. Therefore, in this case, the backwards
output map becomes trivial, leaving only a forwards output map $\beta^o :
\Tt\times X\times A \to B$ and an update map taking inputs in $A$, $\beta^u :
\Tt\times X \times A \to MX$. By filling in the $A$-inputs, we get a system with
$B$-outputs, corresponding to the informal intuition with which we opened this
section: we have a $B$-shaped system with an $A$-shaped hole. Composition of
these systems corresponds to placing systems in parallel using $\otimes$ and
plugging interfaces into holes of the matching shape.

We end this section by briefly sketching the canonical $\otimes$-comonoid
structure on $\HierE$, making $\HierE$ into a `semi-Markov'
\parencite{Fritz2019synthetic} or `copy-discard'
\parencite{Cho2017Disintegration} category. Note that, if a system has the
trivial state space $1$, then (i) tensoring with it is a no-op, and (ii) it has
a trivial update map (assuming that $M1 \cong 1$\footnote{This condition is
satisfied when $M$ is a probability monad like the finite-support distribution
monad, for instance.}). Thus, for each object $Ay$, we obtain a discarding
system $\ground_A : Ay\to 1y$ by taking the trivial state space, trivial update
map, and trivial output map. The copying system $\bcopier_A : Ay\to (A\times
A)y$ again has trivial state space and update map, but now the output map
$\bcopier_A^o : \Tt\times A\to A\times A$ is given by the constant copying map
$(t,a)\mapsto (a,a)$. It is then straightforward to check the comonoid laws.

\section{Dynamical Bayesian inversion}

One consequence of $\HierE$ being a copy-discard category is that we can
instantiate an abstract form of Bayes' rule there, giving rise to a notion of
when one $p\Pa$-coalgebraic system can be seen to be `predicting' or `inverting'
another. In general, Bayes' rule is expressed as an equality between morphisms,
but this is too strong for dynamical systems, which `black-box' their state
spaces: that is to say, we should consider two morphisms (systems) `equal' when
they are observationally equivalent---or, more precisely, when they are related
by a (quasi-)bisimulation.

\begin{defn}
  We define a family of relations \(\sim\) that we collectively call
  \textit{quasi-bisimilarity}.  Given systems \(\vartheta :=
  (X,\vartheta^o,\vartheta^u)\) and \(\psi := (Y,\psi^o,\psi^u)\) in
  \(\pMCoalgT{p}{\Pa}{\Tt}\) and a section \(\sigma\) of \(p\), we first define
  the \textit{trace}\footnote{Note that this is in analogy with the
  \textit{coalgebraic trace}, not the trace of \textit{traced monoidal
    categories}.} or \textit{trajectory} of \(\vartheta\) given \(\sigma\) as
  the morphism
  \[ \mathsf{tr}(\theta,\sigma) := \Tt\times X\xto{\vartheta^o(\mdash)^*\sigma}\sum_{t:\Tt}\sum_{x:X}p[\vartheta^o(t,x)]\xto{(\vartheta^u)^{\rtri_{\Tt}}}\Tt\times\Pa X\xto{\Pa\vartheta^o}\Pa p(1) \, . \]
  Supposing \(\alpha : 1\to\Pa X\) and \(\beta : 1\to\Pa Y\) to be corresponding initial states, we define \(\vartheta \overset{\alpha,\beta}{\sim} \psi\) as the relation
  \[ \vartheta \overset{\alpha,\beta}{\sim} \psi \iff \forall \sigma:\Gamma(p). \, \forall t:\Tt. \, \mathsf{tr}(\vartheta,\sigma)(t)\klcirc\alpha = \mathsf{tr}(\psi,\sigma)(t)\klcirc\beta \, , \]
  where we write \(g\klcirc f\) to indicate Kleisli composition \(g\klcirc f = \mu\circ\Pa g\circ f\) (where \(\mu\) is the multiplication of the monad \(\Pa\)).
  We write \(\vartheta\overset{\exists,\exists}{\sim}\psi\) when there exists some \(\alpha,\beta\) such that \(\vartheta\overset{\alpha,\beta}{\sim}\psi\), and likewise for \(\vartheta\overset{\forall,\forall}{\sim}\psi\), \(\vartheta\overset{\forall,\exists}{\sim}\psi\), and \(\vartheta\overset{\exists,\forall}{\sim}\psi\).
\end{defn}

In light of this definition, we can define an appropriate notion of Bayesian
inversion for $\HierE$:

\begin{defn} \label{def:admit-bayes}
  We say that a system \(c : Xy \to Yy\) in $\HierE$ \textit{admits Bayesian
    inversion} with respect to \(\pi : y \to Xy\), if there exists a system
  \(c^\dag_\pi : Yy \to Xy\) satisfying the equation
  \parencite[eq. 5]{Cho2017Disintegration}:
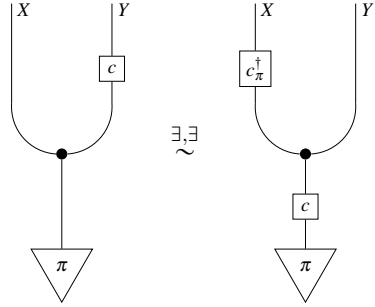
  \begin{equation*} \label{eq:bayes-abstr}
    \scalebox{0.667}{\begin{tikzpicture}
	\begin{pgfonlayer}{nodelayer}
		\node [style=black copier] (1) at (0, 0) {};
		\node [style=none] (2) at (-2, 6) {};
		\node [style=arrow box] (3) at (2, 3.5) {$c$};
		\node [style=none] (4) at (2, 6) {};
		\node [style=none] (5) at (-2, 2) {};
		\node [style=state180] (6) at (0, -4.5) {$\pi$};
		\node [style=none] (7) at (-1.5, 5.75) {$X$};
		\node [style=none] (8) at (2.5, 5.75) {$Y$};
		\node [style=none] (9) at (2, 2) {};
	\end{pgfonlayer}
	\begin{pgfonlayer}{edgelayer}
		\draw [bend left=45] (1) to (5.center);
		\draw (5.center) to (2.center);
		\draw (3) to (4.center);
		\draw (6) to (1);
		\draw [bend right=45] (1) to (9.center);
		\draw (9.center) to (3);
	\end{pgfonlayer}
\end{tikzpicture}
}
    \quad \overset{\exists,\exists}{\sim} \quad
    \scalebox{0.667}{\begin{tikzpicture}
	\begin{pgfonlayer}{nodelayer}
		\node [style=black copier] (0) at (0, 0) {};
		\node [style=none] (1) at (2, 6) {};
		\node [style=arrow box] (2) at (-2, 3.5) {$c^\dag_\pi$};
		\node [style=none] (3) at (-2, 6) {};
		\node [style=none] (4) at (2, 2) {};
		\node [style=state180] (5) at (0, -4.5) {$\pi$};
		\node [style=arrow box] (6) at (0, -2) {$c$};
		\node [style=none] (7) at (-1.5, 5.75) {$X$};
		\node [style=none] (8) at (2.5, 5.75) {$Y$};
		\node [style=none] (9) at (-2, 2) {};
	\end{pgfonlayer}
	\begin{pgfonlayer}{edgelayer}
		\draw [bend right=45] (0) to (4.center);
		\draw (4.center) to (1.center);
		\draw (2) to (3.center);
		\draw (5) to (6);
		\draw (6) to (0);
		\draw [bend left=45] (0) to (9.center);
		\draw (9.center) to (2);
	\end{pgfonlayer}
\end{tikzpicture}
}
  \end{equation*}
  We call \(c^\dag_\pi\) the \textit{Bayesian inversion} of \(c\) with respect
  to \(\pi\), and call the defining relation the \textit{dynamical Bayes' rule}.
\end{defn}

\section{The Laplace doctrine of predictive processing}

In real-world systems, however, even such quasi-bisimulation is too strong.  In
the setting of computational neuroscience, it is proposed
\parencite{Friston2005theory,Buckley2017free} that certain neural circuits
implement \textit{approximate} Bayesian inference by optimizing certain
statistical games \parencite{Smithe2021Compositional1}. A statistical game
consists of a Bayesian lens---a pair of a `forwards' stochastic channel $A\to
\Pa B$ and a `backwards' inversion $\Pa A\times B \to \Pa A$---equipped with a
loss function to evaluate the systems predictive performance. Embodied
predictive systems such as brains then realize these games as dynamical
systems. Here we sketch this functorial semantics, using a category of
`hierarchical bidirectional $\mathsf{Stat}$-systems', following
\parencite{Smithe2020Bayesian,Smithe2021Compositional1}.

We noted above that the category $\Cat{HiBi}$ resembles a category of lenses,
but it does not sufficiently resemble the category of \textit{Bayesian} lenses:
notice that the backwards maps of the latter have codomains of the form $\Pa
A\times T\to \Pa S$ rather than $A\times T\to S$. For this reason, $\Cat{HiBi}$
makes for an inadequate semantic category for predictive processing. However,
all is not lost, for we can define a modification of $\Cat{HiBi}$ by analogy to
the definition of Bayesian lenses as Grothendieck lenses for the indexed
category $\mathsf{Stat}$ of state-dependent maps \parencite{Smithe2020Bayesian}.

\begin{defn}
  Denote by $\Cat{HiBi}_{\Pa}$ the following (semi-)(bi)category. Its objects
  are pairs of objects in $\cat{E}$, and its hom-categories
  $\Cat{HiBi}_{\Pa}((A,S),(B,T))$ are given by $\Cat{Org}^\Tt_{\Pa}(\Pa Ay^S,
  By^T)$. Composition is given by the following family of composite maps:
  \begin{align*}
    & \Cat{HiBi}_{\Pa}((A,S),(B,T)) \times \Cat{HiBi}_{\Pa}((B,T),(C,U)) \\
    & = \Cat{Org}^\Tt_{\Pa}(\Pa Ay^S,By^T) \times \Cat{Org}^\Tt_{\Pa}(\Pa By^T, Cy^U) \\
    & = \pMCoalgT{[\Pa Ay^S,By^T]}{\Pa}{\Tt} \times \pMCoalgT{[\Pa By^T, Cy^U]}{\Pa}{\Tt} \\
    & \to \pMCoalgT{[\Pa Ay^S,By^T]\otimes[\Pa By^T, Cy^U]}{\Pa}{\Tt} \\
    & \to \pMCoalgT{[\Pa Ay^S,\Pa By^T]\otimes[\Pa By^T, Cy^U]}{\Pa}{\Tt} \\
    & \to \pMCoalgT{[\Pa Ay^S, Cy^U]}{\Pa}{\Tt} \\
    & = \Cat{HiBi}_{\Pa}((A,S),(C,U))
  \end{align*}
  where the fourth line is generated from the monadic unit $\eta_B : B\to \Pa B$
  by $\Coalg([\Pa y^S,(\eta_B)y^T])$.
\end{defn}

\begin{rmk}
  Note that we say 'semi-'(bi)category: this is because $\Cat{HiBi}_{\Pa}$ does
  not have identities. This is not problematic for our work here; and of course
  $\Cat{Org}^\Tt_{\Pa}$ itself does have identities.
\end{rmk}

We are now in a position to sketch the `Laplace doctrine' of dynamical semantics
for approximate inference. We first recall the notion of $D$-Bayesian inference
game \parencite{Smithe2021Compositional1}:

\begin{defn}[Bayesian inference] \label{def:bayes-game}
  Let \(D : \Kl(\Pa)(I,X) \times \Kl(\Pa)(1,X) \to \rr\) be a measure of
  divergence between states on \(X\). Then a (simple) \(D\)\textit{-Bayesian
    inference} game is a statistical game \((X,X) \to (Y,Y)\) with fitness
  function \(\phi : \Kl(\Pa)(1, X)\times \Kl(\Pa)(Y, X) \to \rr\) given by
  \[\phi(\pi,k) = \E_{y\sim k\klcirc c\klcirc\pi} \left[ D\left(c'_\pi(y),
    c^\dag_\pi(y)\right) \right]\] where \((c,c')\) constitutes the lens part
  of the game and \(c^\dag_\pi\) is the exact inversion of \(c\) with respect to
  \(\pi\).
\end{defn}

Write $D_{KL}$ for the Kullback-Leibler divergence. Given a $D_{KL}$-Bayesian
inference game $(\gamma,\rho,\phi) : (X,X)\to(Y,Y)$ where $X$ and $Y$ are
Euclidean spaces and whose forward and backward channels are constrained to
output Gaussian distributions, the Laplace doctrine returns a hierarchical
bidirectional $\mathsf{Stat}$-system minimizing an upper bound on the divergence
between each approximate posterior $\rho_\pi$ and the `true' posterior
$\gamma^\dag_\pi$, for any Gaussian state $\pi : \Pa X$.

\begin{rmk}
  Note that the statistical properties of the system are not the focus of this
  paper: this doctrine is merely being used to illustrate the coalgebraic framework.
\end{rmk}

The Laplace doctrine hinges on the following approximation, whose proof we defer
to \ref{proof:laplace-approx}.

\begin{lemma}[Laplace approximation] \label{lemma:laplace-approx}
  Given a $D_{KL}$-Bayesian inference game $(\gamma,\rho,\phi):(X,X)\to(Y,Y)$
  with forwards channel $\gamma : X\to \Pa Y$ constrained to emit Gaussian
  distributions, write $\mu_\gamma(x) : \rr^{|Y|}$ for the mean of $\gamma(x)$
  and $\Sigma_\gamma(x) : \rr^{|Y|\times |Y|}$ for its covariance matrix, and
  assume that for all \(y : Y\), the eigenvalues of \(\Sigma_{\rho_\pi}(y)\) are
  small.

  Then the loss \(\phi : \Kl(\Pa)(1, X)\times \Kl(\Pa)(Y, X) \to \rr\) is
  approximately bounded from above by
  \begin{align*}
  \phi(\pi, k)
  &= \E_{y\sim k\klcirc \gamma\klcirc\pi} \left[ D\left(\rho_\pi(y), \gamma^\dag_\pi(y)\right) \right] \\
  &\leq \E_{y\sim k\klcirc \gamma\klcirc\pi} \left[ D\left(\rho_\pi(y), \gamma^\dag_\pi(y)\right) - \log p_{\gamma\klcirc\pi}(y) \right] \\
  &= \E_{y\sim k\klcirc \gamma\klcirc\pi} \big[ \Fa(y) \big]
   \approx \E_{y \sim k \klcirc \gamma \klcirc \pi} \big[ \Fa^L(y) \big]
  \end{align*}
  where $\Fa$ is called the \textit{free energy} and where $\Fa^L$ is its
  \textit{Laplace approximation},
  \begin{align} \label{eq:laplace-energy}
    \Fa^L(y)
    & = E_{(\pi,\gamma)}\left(\mu_{\rho_\pi}(y), y\right) - S_X \left[ \rho_\pi(y) \right] \\
    & = -\log p_\gamma(y|\mu_{\rho_\pi}(y)) -\log p_\pi(\mu_{\rho_\pi}(y)) - S_X \left[ \rho_\pi(y) \right] \nonumber
  \end{align}
  where \(S_x[\rho_\pi(y)] = \E_{x \sim \rho_\pi(y)} [ -\log p_{\rho_\pi}(x|y)
  ]\) is the Shannon entropy of \(\rho_\pi(y)\), and \(p_\gamma : Y \times X \to
  [0,1]\), \(p_\pi : X \to [0,1]\), and \(p_{\rho_\pi} : X \times Y \to [0,1]\)
  are density functions for \(\gamma\), \(\pi\), and \(\rho_\pi\) respectively.
  The approximation is valid when \(\Sigma_{\rho_\pi}\) satisfies
  \begin{equation} \label{eq:laplace-sigma-rho-pi}
    \Sigma_{\rho_\pi} (y) = \left(\partial_x^2 E_{(\pi,\gamma)}\right)\left( \mu_{\rho_\pi}(y), y\right)^{-1} \, .
  \end{equation}
\end{lemma}

With this approximation in hand, and given such a statistical game \((\gamma,
\rho, \phi)\), we will construct a hierarchical bidirectional
$\mathsf{Stat}$-system \(\Fun{Laplace}(\gamma, \rho, \phi)\) performing
approximate stochastic gradient descent on the loss function, with respect to
the statistical parameters of the inversions \(\rho_\pi\). We will work in
discrete time, $\Tt = \nn$, although all of what follows can be done in
continuous time, $\Tt = \rr_+$, by replacing the discrete update steps by
stochastic differential equations.

Since the entropy \(S_X \left[ \rho_\pi(y) \right]\) depends only on the variance \(\Sigma_{\rho_\pi} (y)\), to optimize the mean \(\mu_{\rho_\pi} (y)\) it suffices to consider only the energy \(E_{(\pi,\gamma)}(\mu_{\rho_\pi}(y), y)\). We have
\begin{align*}
E_{(\pi,\gamma)}(x, y)
&= - \log p_\gamma(y|x) - \log p_\pi(x) \\
&= - \frac{1}{2} \innerprod{\epsilon_\gamma(y,x)}{{\Sigma_\gamma(x)}^{-1} {\epsilon_\gamma(y,x)}}
   - \frac{1}{2} \innerprod{\epsilon_\pi(x)}{{\Sigma_\pi}^{-1} {\epsilon_\pi(x)}} \\
&\quad
   + \log \sqrt{(2 \pi)^{|Y|} \det \Sigma_\gamma(x) }
   + \log \sqrt{(2 \pi)^{|X|} \det \Sigma_\pi }
\end{align*}
and a straightforward computation shows that
\[
\partial_{x} E_{(\pi,\gamma)}(x, y)
= - \partial_x \mu_\gamma (x)^T {\Sigma_\gamma(x)}^{-1} \epsilon_\gamma(y,x) + {\Sigma_\pi}^{-1} \epsilon_\pi(x) \, .
\]
Let \(\eta_\gamma(y, x) := {\Sigma_\gamma(x)}^{-1} \epsilon_\gamma(y, x)\) and \(\eta_\pi(x) := {\Sigma_\pi}^{-1} \epsilon_\pi(x)\), so that
\begin{equation} \label{eq:laplace-dEdx}
\partial_{x} E_{(\pi,\gamma)}(x, y)
= - \partial_x \mu_\gamma (x)^T \eta_\gamma(y,x) + \eta_\pi(x) \, .
\end{equation}
Note that \(E_{(\pi,\gamma)}\) defines a function \(X \times Y \to \rr\).
We will use the domain $X\times Y$ of this function as the state space of our
system. To avoid ambiguity, we will write $\overrightarrow{X}$ to indicate the
space $X$ when it is used as an input in the `forwards' direction, and
$\overleftarrow{Y}$ to indicate the space $Y$ when it is used as an input in the
`backwards' direction.

Our system \(\Fun{Laplace}(\gamma, \rho, \phi)\) will therefore have the type
\begin{align*}
  (& X\times Y, \quad \beta^o_1 : X\times Y\times \Pa\overrightarrow{X} \to \overrightarrow{Y}, \\
  & \beta^o_2 : X\times Y\times \Pa\overrightarrow{X} \times \overleftarrow{Y} \to \overleftarrow{X}, \\
  & \beta^u : X\times Y\times \Pa\overrightarrow{X} \times \overleftarrow{Y} \to \Pa(X\times Y)).
\end{align*}
We define $\beta^o_1$ to be the projection of the second factor $Y$ of the state
space onto $Y$, and $\beta^o_2$ to be the projection of the first factor $X$
onto $X$. The update map $\beta^u : X\times Y\times \Pa\overrightarrow{X} \times
\overleftarrow{Y} \to \Pa(X\times Y)$ is then given by composing the
commutativity (or `double strength') of the monad $\Pa$, $\mathsf{dst} : \Pa
X\times \Pa Y\to\Pa(X\times Y)$, after the following map (represented as a
string diagram in $\cat{E}$):
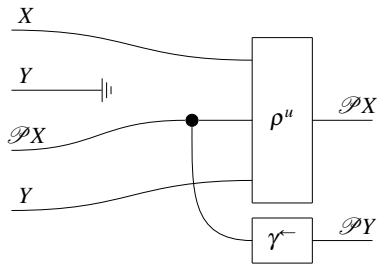
\[
\scalebox{0.8}{\begin{tikzpicture}
	\begin{pgfonlayer}{nodelayer}
		\node [style=none] (0) at (2, 3.75) {};
		\node [style=none] (1) at (4, 3.75) {};
		\node [style=none] (2) at (2, -1.75) {};
		\node [style=none] (3) at (4, -1.75) {};
		\node [style=none] (4) at (2, 3) {};
		\node [style=black copier] (5) at (0, 1) {};
		\node [style=none] (6) at (2, -1) {};
		\node [style=none] (7) at (4, 1) {};
		\node [style=none] (13) at (3, 1) {$\rho^u$};
		\node [style=none] (16) at (2, -2.25) {};
		\node [style=none] (17) at (4, -2.25) {};
		\node [style=none] (18) at (4, -3.75) {};
		\node [style=none] (19) at (2, -3.75) {};
		\node [style=none] (20) at (3, -3) {$\gamma^\leftarrow$};
		\node [style=none] (21) at (2, -3) {};
		\node [style=none] (22) at (4, -3) {};
		\node [style=none] (23) at (6, 1) {};
		\node [style=none] (24) at (6, -3) {};
		\node [style=none] (45) at (-6, 4) {};
		\node [style=none] (46) at (-6, 0) {};
		\node [style=none] (47) at (-6, -2) {};
		\node [style=none] (50) at (5.5, 1.5) {$\mathcal{P} X$};
		\node [style=none] (51) at (5.5, -2.5) {$\mathcal{P} Y$};
		\node [style=none] (52) at (-5.5, -1.5) {$Y$};
		\node [style=none] (53) at (-5.5, 0.5) {$\Pa X$};
		\node [style=none] (54) at (-5.5, 4.5) {$X$};
		\node [style=none] (56) at (2, 1) {};
		\node [style=none] (57) at (-6, 2) {};
		\node [style=none] (58) at (-5.5, 2.5) {$Y$};
		\node [style=discarder0] (59) at (-3, 2) {};
	\end{pgfonlayer}
	\begin{pgfonlayer}{edgelayer}
		\draw (0.center) to (1.center);
		\draw (1.center) to (3.center);
		\draw (3.center) to (2.center);
		\draw (2.center) to (0.center);
		\draw (16.center) to (17.center);
		\draw (17.center) to (18.center);
		\draw (18.center) to (19.center);
		\draw (19.center) to (16.center);
		\draw (7.center) to (23.center);
		\draw (22.center) to (24.center);
		\draw [in=180, out=0, looseness=1.25] (47.center) to (6.center);
		\draw (5) to (56.center);
		\draw [in=180, out=-90, looseness=1.25] (5) to (21.center);
		\draw [in=-180, out=0, looseness=1.25] (46.center) to (5);
		\draw (57.center) to (59);
		\draw [in=180, out=0, looseness=1.25] (45.center) to (4.center);
	\end{pgfonlayer}
\end{tikzpicture}
}
\]
where \((-)^\leftarrow := \mu^{\Pa} \circ \Pa(-)\) denotes Kleisli extension
(for \(\mu^{\Pa}\) the multiplication of the monad \(\Pa\)), so that
\(\gamma^\leftarrow := \mu^{\Pa}_Y \circ \Pa(\gamma) : \Pa X \to \Pa Y\).

In turn, the map \(\rho^u : X \times \Pa X \times Y \to \Pa X\) is defined by
\begin{align*}
\rho^u : X \times \Pa X \times Y & \to \rr^{|X|} \times \rr^{|X|\times|X|}
\hookrightarrow \Pa X \\ (x, \pi, y) & \mapsto \left(x - \lambda \partial_{x}
E_{(\pi,\gamma)}(x, y), \Sigma_{\rho_\pi}^* (y)\right)
\end{align*}
where the inclusion into \(\Pa X\) picks the Gaussian state with the given
statistical parameters, where \(\lambda : \rr_+\) is some choice of ``learning
rate'', where \(\Sigma_{\rho_\pi}^* (y)\) is as above and in Equation
\eqref{eq:laplace-sigma-rho-pi}, and where \(\partial_{x} E_{(\pi,\gamma)}(x,
y)\) is as in Equation \eqref{eq:laplace-dEdx}.

Observe that the factor \(\rho^u\) performs approximate stochastic gradient
descent on the free energy: for a given input \(y:Y\), the mean trajectory of
the system follows the update law \(\mu_{\rho} \mapsto \mu_{\rho} - \lambda
\partial_{\mu_{\rho}} E_{(\pi,\gamma)}(\mu_{\rho}, y)\), and, when
\(\Sigma_{\rho} (y) = \Sigma_{\rho}^* (y)\), we have \(\partial_{\mu_{\rho}}
E_{(\pi,\gamma)}(\mu_{\rho}, y) \approx \partial_{\mu_{\rho}} \Fa(y)\). Note
also that the update map \(\rho^u\) depends on a prior, just as the inversion
map \(\rho\) of the lens \((\gamma, \rho)\) does.

A full treatment of the Laplace doctrine will appear in a forthcoming sequeal to
the author's \parencite{Smithe2021Compositional1}.

\section{Conclusions; current and future work}

In this work we have sketched a framework for treating open dynamical systems of
a general nature as coalgebras for certain polynomial functors or---in the case
of systems with side-effects such as randomness---certain generalizations
thereof. Although we have attempted to give a wide overview of the applicability
of these structures, with a particular focus on the adaptive systems of primary
interest to the author, we are aware that we have barely scratched the surface
of their use and relationships. Here, we briefly list some avenues of current
and future work.

Our current principal focus is on exploring the connections between these
structures and other compositional treatments of dynamical systems.  In
particular, relating our categories to the respective frameworks of Myers
\parencite{Myers2020Double}, Libkind \parencite{Libkind2020Algebra} and Baez and
colleagues (\textit{e.g.}, \parencite{Baez2021Structured}).  Evidently, the
structures presented here are most closely in line with the approaches explored
by Spivak \parencite{Spivak2021Learners,Schultz2019Dynamical}, and are
particularly interested in generalizing his topos-theoretic perspective: given
that the category of discrete-time deterministic systems over a polynomial $p$
forms a topos, we suspect that so too does $\Coalg^{\Tt}(p)$. We are also
seeking the connections between these putative topoi and the topoi of behaviour
types \parencite{Schultz2019Dynamical} as well as with coalgebraic logic
\parencite{Jacobs2017Introduction}, particularly in its modal forms. We hope
that we can further develop the theory of $\PolyM{M}$ to support some of these
methods, too.

Finally, there are a number of ways in which this framework should be made more
elegant. In particular, we hope to cast a number of properties instead as
structures, including the comonoid-homomorphism property of our main definition,
and the explicit definitions of random and bundle dynamical systems. With
particular respect to the latter, we expect there to be an inductive story of
nested parameterization, which appears to the author to have an opetopic shape
closely connected to the $\Cat{Para}$ construction
\parencite{Capucci2021Parameterized}.

\section{References}
{\printbibliography[heading=none]}

\appendix
\section{Extra proofs and structures}

\begin{prop} \label{prop:poly-dyn-idx}
  \(\Coalg^\Tt(p)\) extends to a polynomially-indexed category, \(\Coalg^\Tt :
  \Poly{E} \to \Cat{Cat}\). Suppose \(\varphi : p \to q\) is a morphism of
  polynomials.  We define a corresponding functor \(\Coalg^\Tt(\varphi) :
  \Coalg^\Tt(p) \to \Coalg^\Tt(q)\) as follows.  Suppose \((X, \vartheta^o,
  \vartheta^u) : \Coalg^\Tt(p)\) is an object (dynamical system) in
  \(\Coalg^\Tt(p)\). Then \(\Coalg^\Tt(\varphi)(X, \vartheta^o, \vartheta^u)\) is
  defined as the triple \((X, \varphi_1 \circ \vartheta^o, \vartheta^u \circ
  {\vartheta^o}^\ast \varphi^\#) : \Coalg^\Tt(q)\), where the two maps are
  explicitly the following composites:
  \begin{gather*}
    \Tt \times X \xto{\vartheta^o} p(1) \xto{\varphi_1} q(1) \, ,
    \qquad
    \Sum_{t:\Tt} \Sum_{x:X} q[\varphi_1 \circ \vartheta^o(t, x)] \xto{{\vartheta^o}^\ast \varphi^\#} \Sum_{t:\Tt} \Sum_{x:X} p[\vartheta^o(t, x)] \xto{\vartheta^u} X \, .
  \end{gather*}
  On morphisms, \(\Coalg^\Tt(\varphi)(f) : \Coalg^\Tt(\varphi)(X, \vartheta^o,
  \vartheta^u) \to \Coalg^\Tt(\varphi)(Y, \psi^o, \psi^u)\) is given by the same
  underlying map \(f : X \to Y\) of state spaces.
  \begin{proof}
    We need to check that \(\Coalg^\Tt(\varphi)(X, \vartheta^o, \vartheta^u)\)
    satisfies the flow conditions of Definition \ref{def:poly-dyn}, that
    \(\Coalg^\Tt(\varphi)(f)\) satisfies the naturality condition of Proposition
    \ref{prop:poly-dyn}, and that \(\Coalg^\Tt\) is functorial with respect to
    polynomials. We begin with the flow condition. Given a section \(\tau : q(1)
    \to \Sum_{j:q(1)} q[j]\) of \(q\), we require the closures
    \(\Coalg^\Tt(\varphi)(\vartheta)^\tau : \Tt \times X \to X\) given by
    \[
    \Sum_{t:\Tt} X \xto{\vartheta^o(-)^\ast \tau} \Sum_{t:\Tt} \Sum_{x:X} q[\varphi_1 \circ \vartheta^o(t, x)] \xto{{\vartheta^o}^\ast \varphi^\#} \Sum_{t:\Tt} \Sum_{x:X} p[\vartheta^o(t, x)] \xto{\vartheta^u} X
    \]
    to satisfy \(\Coalg^\Tt(\varphi)(\vartheta)^\tau (0) = \id_X\) and
    \(\Coalg^\Tt(\varphi)(\vartheta)^\tau (s+t) =
    \Coalg^\Tt(\varphi)(\vartheta)^\tau (s) \circ
    \Coalg^\Tt(\varphi)(\vartheta)^\tau (t)\). Note that the following diagram
    commutes, by the definition of \(\varphi^\#\),
    \[\begin{tikzcd}
	{\Sum_{i:p(1)} p[i]} && {\Sum_{i:p(1)}q[\varphi_1(i)]} && {p(1)} \\
	\\
	{p(1)} && {p(1)}
	\arrow["{\varphi_1^\ast q}"', from=1-3, to=3-3]
	\arrow[Rightarrow, no head, from=3-3, to=1-5]
	\arrow["{\varphi_1^\ast \tau}"', from=1-5, to=1-3]
	\arrow["{\varphi^\#}"', from=1-3, to=1-1]
	\arrow[Rightarrow, no head, from=3-1, to=3-3]
	\arrow["p"', from=1-1, to=3-1]
    \end{tikzcd}\]
    so that \(\varphi^\# \circ \varphi_1^\ast \tau\) is a section of
    \(p\). Therefore, letting \(\sigma := \varphi^\# \circ \varphi_1^\ast
    \tau\), for \(\Coalg^\Tt(\varphi)(\vartheta)^\tau\) to satisfy the flow
    condition for \(\tau\) reduces to \(\vartheta^\sigma\) satisfying the flow
    condition for \(\sigma\). But this is given \textit{ex hypothesi} by
    Definition \ref{def:poly-dyn}, for any such section \(\sigma\), so
    \(\Coalg^\Tt(\varphi)(\vartheta)^\tau\) satisfies the flow condition for
    \(\tau\). And since \(\tau\) was any section, we see that
    \(\Coalg^\Tt(\varphi)(\vartheta)\) satisfies the flow condition generally.

    The proof that \(\Coalg^\Tt(\varphi)(f)\) satisfies the naturality condition
    of Proposition \ref{prop:poly-dyn} proceeds similarly.  Supposing again that
    \(\tau\) is any section of \(q\), we require the following diagram to
    commute for any time \(t : \Tt\):
    \[\begin{tikzcd}
	X && {\Sum_{x:X} q[\varphi_1 \circ \vartheta^o(t, x)]} && {\Sum_{x:X} p[\vartheta^o(t, x)]} && X \\
	\\
	Y && {\Sum_{y:Y} q[\varphi_1 \circ \psi^o(t, x)]} && {\Sum_{y:Y} p[\psi^o(t, x)]} && Y
	\arrow["f"', from=1-1, to=3-1]
	\arrow["f", from=1-7, to=3-7]
	\arrow["{\vartheta^o(t)^\ast \varphi_1^\ast \tau}", from=1-1, to=1-3]
	\arrow["{\vartheta^o(t)^\ast \varphi^\#}", from=1-3, to=1-5]
	\arrow["{\vartheta^u(t)}", from=1-5, to=1-7]
	\arrow["{\psi^o(t)^\ast \varphi_1^\ast \tau}", from=3-1, to=3-3]
	\arrow["{\vartheta^o(t)^\ast \varphi^\#}", from=3-3, to=3-5]
	\arrow["{\psi^u(t)}", from=3-5, to=3-7]
    \end{tikzcd}\]
    Again letting \(\sigma := \varphi^\# \circ \varphi_1^\ast \tau\), we see
    that this diagram reduces exactly to the diagram in Proposition
    \ref{prop:poly-dyn} by the functoriality of pullback, and since \(f\) makes
    that diagram commute, it must also make this diagram commute.

    Finally, to show that \(\Coalg^\Tt\) is functorial with respect to polynomials
    amounts to checking that composition and pullback are functorial; but this
    is a basic result of category theory.
  \end{proof}
\end{prop}

\begin{prop} \label{prop:dyn-n-coalg}
  When \(\Tt = \nn\), the category \(\Coalg^\nn(p)\) of open dynamical systems
  over \(p\) with time \(\nn\) is equivalent to the topos \(p\mdash\Cat{Coalg}\)
  of \(p\)-coalgebras \parencite{Spivak2021Learners}.
  \begin{proof}
    \(p\mdash\Cat{Coalg}\) has as objects pairs \((S, \beta)\) where \(S :
    \cat{E}\) is an object in \(\cat{E}\), \(\beta : S \to p \triangleleft S\)
    is a morphism of polynomials (interpreting \(S\) as the constant copresheaf
    on the set \(S\)), and \(\triangleleft\) denotes the composition monoidal
    product in \(\Poly{E}\) (\textit{i.e.}, composing the corresponding
    copresheaves \(\cat{E} \to \cat{E}\)). A straightforward computation shows
    that, interpreted as an object in \(\cat{E}\), \(p \triangleleft S\)
    corresponds to \(\sum_{i:p(1)} S^{\,p[i]}\).  By the universal property of
    the dependent sum, a morphism \(\beta : S \to \sum_{i:p(1)} S^{\,p[i]}\)
    therefore corresponds bijectively to a pair of maps \(\beta^o : S \to p(1)\)
    and \(\beta^u : \sum_{s:S} p[\beta^o(s)] \to X\).  By Proposition
    \ref{prop:open-transition-map}, such a pair is equivalently a discrete-time
    open dynamical system over \(p\) with state space \(S\): that is, the
    objects of \(p\mdash\Cat{Coalg}\) are in bijection with those of
    \(\Coalg^\nn(p)\).

    Next, we show that the hom-sets \(p\mdash\Cat{Coalg}\big((S, \beta), (S',
    \beta')\big)\) and \(\Coalg^\nn(p)\big((S, \beta^o, \beta^u), (S', \beta'^o,
    \beta'^u)\big)\) are in bijection.  A morphism \(f : (S, \beta) \to (S',
    \beta')\) of \(p\)-coalgebras is a morphism \(f : S \to S'\) between the
    state spaces such that \(\beta' \circ f = (p \triangleleft f) \circ
    \beta\). Unpacking this, we find that this means the following diagram in
    \(\cat{E}\) must commute for any section \(\sigma\) of \(p\):
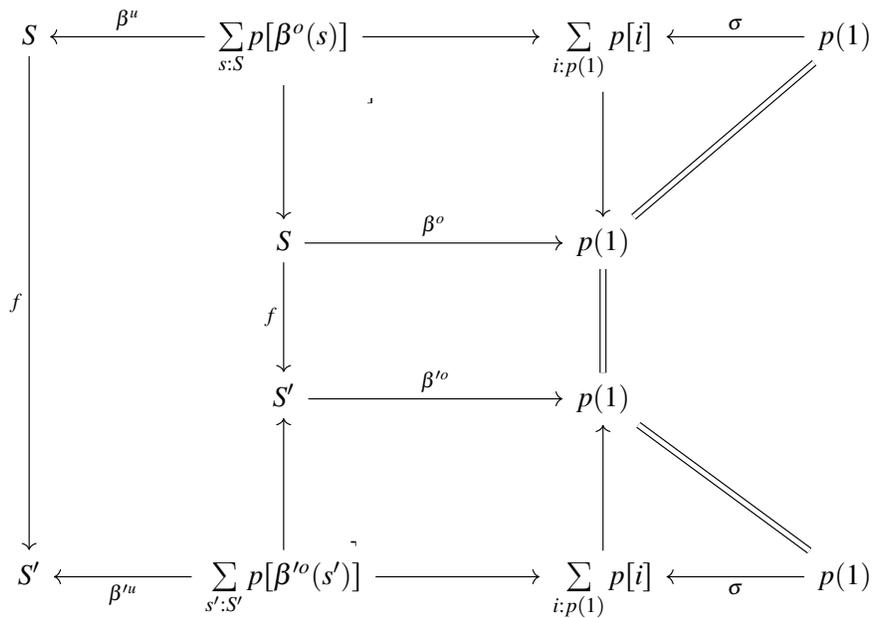
    \[\begin{tikzcd}
	S && {\Sum_{s:S} p[\beta^o(s)]} && {\Sum_{i:p(1)} p[i]} && {p(1)} \\
	&&& {} \\
	&& S && {p(1)} \\
	\\
	&& {S'} && {p(1)} \\
	&&& {} \\
	{S'} && {\Sum_{s':S'} p[\beta'^o(s')]} && {\Sum_{i:p(1)} p[i]} && {p(1)}
	\arrow[Rightarrow, no head, from=3-5, to=5-5]
	\arrow["f"', from=1-1, to=7-1]
	\arrow["f"', from=3-3, to=5-3]
	\arrow[from=1-3, to=1-5]
	\arrow["{\beta^o}", from=3-3, to=3-5]
	\arrow[from=1-5, to=3-5]
	\arrow[from=1-3, to=3-3]
	\arrow["\lrcorner"{anchor=center, pos=0.125}, draw=none, from=1-3, to=2-4]
	\arrow[from=7-3, to=7-5]
	\arrow[from=7-5, to=5-5]
	\arrow["{\beta'^o}", from=5-3, to=5-5]
	\arrow[from=7-3, to=5-3]
	\arrow["\lrcorner"{anchor=center, pos=0.125, rotate=90}, draw=none, from=7-3, to=6-4]
	\arrow["{\beta^u}"', from=1-3, to=1-1]
	\arrow["{\beta'^u}", from=7-3, to=7-1]
	\arrow[Rightarrow, no head, from=3-5, to=1-7]
	\arrow[Rightarrow, no head, from=5-5, to=7-7]
	\arrow["\sigma", from=7-7, to=7-5]
	\arrow["\sigma"', from=1-7, to=1-5]
    \end{tikzcd}\]
    Pulling the arbitrary section \(\sigma\) back along the `output' maps
    \(\beta^o\) and \(\beta'^o\) means that the following commutes:
    
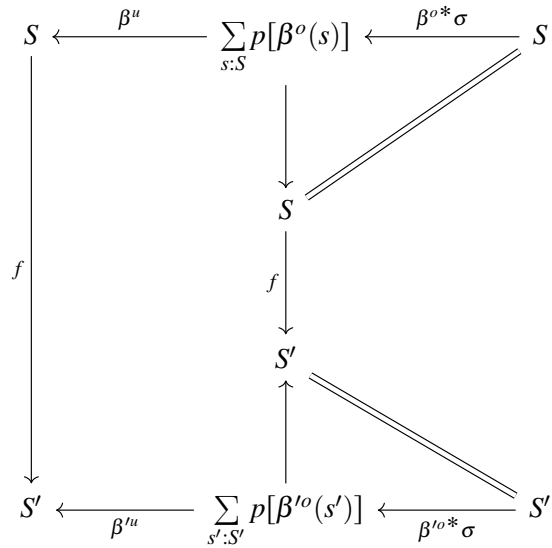
\[\begin{tikzcd}
	S && {\Sum_{s:S} p[\beta^o(s)]} && S \\
	\\
	&& S \\
	\\
	&& {S'} \\
	\\
	{S'} && {\Sum_{s':S'} p[\beta'^o(s')]} && {S'}
	\arrow["f"', from=1-1, to=7-1]
	\arrow["f"', from=3-3, to=5-3]
	\arrow[from=1-3, to=3-3]
	\arrow[from=7-3, to=5-3]
	\arrow["{\beta^u}"', from=1-3, to=1-1]
	\arrow["{\beta'^u}", from=7-3, to=7-1]
	\arrow[Rightarrow, no head, from=3-3, to=1-5]
	\arrow[Rightarrow, no head, from=5-3, to=7-5]
	\arrow["{{\beta^o}^\ast \sigma}"', from=1-5, to=1-3]
	\arrow["{{\beta'^o}^\ast \sigma}", from=7-5, to=7-3]
    \end{tikzcd}\]
    Forgetting the vertical projections out of the pullbacks gives:
    
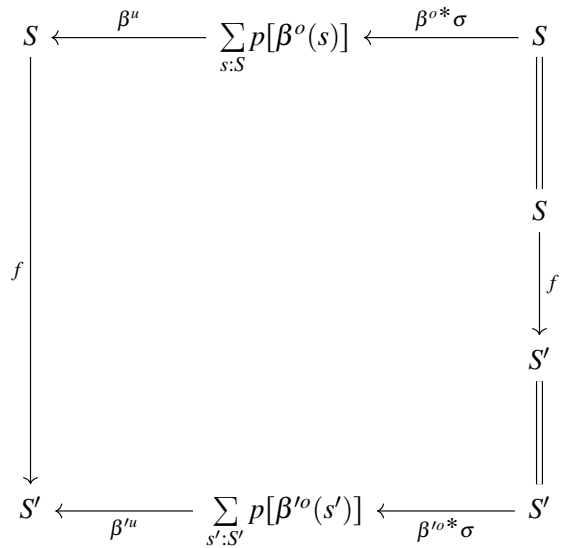
\[\begin{tikzcd}
	S && {\Sum_{s:S} p[\beta^o(s)]} && S \\
	\\
	&&&& S \\
	\\
	&&&& {S'} \\
	\\
	{S'} && {\Sum_{s':S'} p[\beta'^o(s')]} && {S'}
	\arrow["f"', from=1-1, to=7-1]
	\arrow["{\beta^u}"', from=1-3, to=1-1]
	\arrow["{\beta'^u}", from=7-3, to=7-1]
	\arrow["{{\beta^o}^\ast \sigma}"', from=1-5, to=1-3]
	\arrow["{{\beta'^o}^\ast \sigma}", from=7-5, to=7-3]
	\arrow[Rightarrow, no head, from=5-5, to=7-5]
	\arrow["f", from=3-5, to=5-5]
	\arrow[Rightarrow, no head, from=3-5, to=1-5]
    \end{tikzcd}\]
    Finally, by collapsing the identity maps and reflecting the diagram
    horizontally, we obtain
    
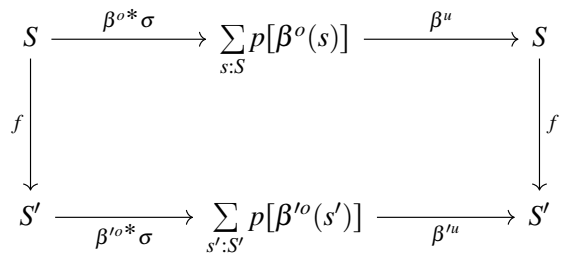
\[\begin{tikzcd}
	S && {\Sum_{s:S} p[\beta^o(s)]} && S \\
	\\
	{S'} && {\Sum_{s':S'} p[\beta'^o(s')]} && {S'}
	\arrow["f"', from=1-1, to=3-1]
	\arrow["f", from=1-5, to=3-5]
	\arrow["{{\beta^o}^\ast \sigma}", from=1-1, to=1-3]
	\arrow["{{\beta'^o}^\ast \sigma}"', from=3-1, to=3-3]
	\arrow["{\beta^u}", from=1-3, to=1-5]
	\arrow["{\beta'^u}"', from=3-3, to=3-5]
    \end{tikzcd}\]
    which we recognize from Proposition \ref{prop:poly-dyn} as the defining
    characteristic of a morphism in \(\Coalg^\nn(p)\). Finally, we note that each
    of these steps is bijective, and so we have the desired bijection of
    hom-sets.
  \end{proof}
\end{prop}

\begin{defn}[Category of open random dynamical systems over \(p\)] \label{def:cat-open-rds}
  Writing \(\vartheta := (\pi_X, \vartheta^o, \vartheta^u)\) and
  \(\psi := (\pi_Y, \psi^o, \psi^u)\), a morphism \(f : \vartheta \to \psi\) is
  a map \(f: X \to Y\) in \(\cat{E}\) making the following diagram commute for
  all times \(t : \Tt\) and sections \(\sigma\) of \(p\):
  \[\begin{tikzcd}
	X &&&& {\Sum_{x:X} p[\vartheta^o(t, x)]} &&&& X \\
	\\
	&& \Omega &&&& \Omega \\
	\\
	Y &&&& {\Sum_{y:Y} p[\psi^o(t, y)]} &&&& Y
	\arrow["{\pi_X}"', from=1-1, to=3-3]
	\arrow["{\pi_X}", from=1-9, to=3-7]
	\arrow["{\theta(t)}", from=3-3, to=3-7]
	\arrow["{\vartheta^o(t)^\ast \sigma}", from=1-1, to=1-5]
	\arrow["{\vartheta^u(t)}", from=1-5, to=1-9]
	\arrow["{\psi^o(t)^\ast \sigma}"', from=5-1, to=5-5]
	\arrow["{\psi^u(t)}"', from=5-5, to=5-9]
	\arrow["{\pi_Y}", from=5-1, to=3-3]
	\arrow["{\pi_Y}"', from=5-9, to=3-7]
	\arrow["f"', from=1-1, to=5-1]
	\arrow["f", from=1-9, to=5-9]
  \end{tikzcd}\]
  Identities are given by the identity maps on state-spaces. Composition is
  given by pasting of diagrams.
\end{defn}

\begin{prop}[Opindexed category of open random dynamical systems over polynomials] \label{prop:poly-rdyn-idx-details}
  By the universal property of the product \(\times\) in \(\Cat{Cat}\), it
  suffices to define the actions of \(\RDynT{T}\) separately on morphisms of
  polynomials and on morphisms of closed measure-preserving systems.

  Suppose therefore that \(\varphi : p \to q\) is a morphism of
  polynomials. Then, for each measure-preserving system \((\theta, \beta) :
  \Cat{Cat}(\deloop{\Tt}, \cat{E})_{\Pa}\), we define the functor
  \(\RDynT{T}(\varphi, \theta) : \RDynT{T}(p, \theta) \to \RDynT{T}(q, \theta)\)
  as follows. Let \(\vartheta := (\pi_X : X \to \Omega, \vartheta^o,
  \vartheta^u) : \RDynT{T}(p, \theta)\) be an object (open random dynamical
  system) in \(\RDynT{T}(p, \theta)\). Then \(\RDynT{T}(\varphi,
  \theta)(\vartheta)\) is defined as the triple \((\pi_X, \varphi_1 \circ
  \vartheta^o, \vartheta^u \circ {\varphi^o}^\ast \varphi^\#) : \RDynT{T}(q,
  \theta)\), where the two maps are explicitly the following composites:
  \begin{gather*}
    \Tt \times X \xto{\vartheta^o} p(1) \xto{\varphi_1} q(1) \, ,
    \qquad
    \Sum_{t:\Tt} \Sum_{x:X} q[\varphi_1 \circ \vartheta^o(t, x)] \xto{{\vartheta^o}^\ast \varphi^\#} \Sum_{t:\Tt} \Sum_{x:X} p[\vartheta^o(t, x)] \xto{\vartheta^u} X \, .
  \end{gather*}
  On morphisms \(f : (\pi_X : X \to \Omega, \vartheta^o, \vartheta^u) \to (\pi_Y
  : Y \to \Omega, \psi^o, \psi^u)\), the image \[\RDynT{T}(\varphi, \theta)(f) :
  \RDynT{T}(\varphi, \theta)(\pi_X, \vartheta^o, \vartheta^u) \to
  \RDynT{T}(\varphi, \theta)(\pi_Y, \psi^o, \psi^u)\] is given by the same
  underlying map \(f : X \to Y\) of state spaces.

  Next, suppose that \(\phi : (\theta, \beta) \to (\theta', \beta')\) is a
  morphism of closed measure-preserving dynamical systems, and let \(\Omega' :=
  \theta'(\ast)\) be the state space of the system \(\theta'\). By Proposition
  \ref{prop:metric-sys}, the morphism \(\phi\) corresponds to a map \(\phi :
  \Omega \to \Omega'\) on the state spaces that preserves both flow and
  measure. Therefore, for each polynomial \(p : \Poly{E}\), we define the
  functor \(\RDynT{T}(p, \phi) : \RDynT{T}(p, \theta) \to \RDynT{T}(p,
  \theta')\) by post-composition. That is, suppose given open random dynamical
  systems and morphisms over \((p, \theta)\) as in the diagram of Proposition
  \ref{prop:poly-rdyn}. Then \(\RDynT{T}(p, \phi)\) returns the following
  diagram:
  
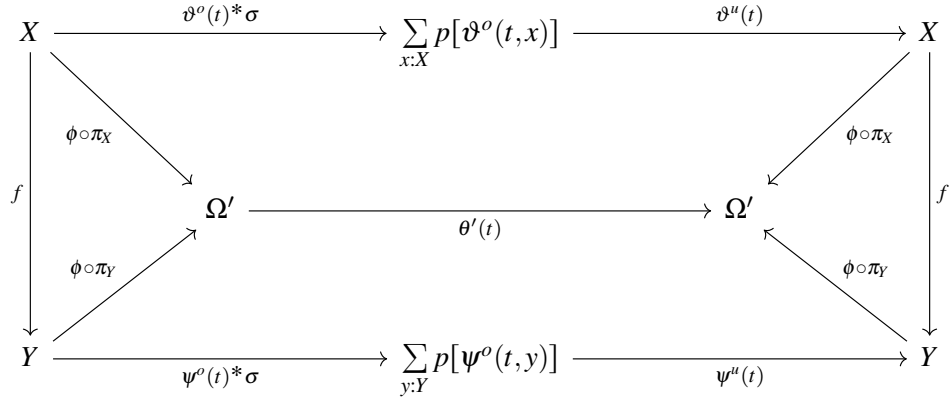
\[\begin{tikzcd}
	X &&&& {\Sum_{x:X} p[\vartheta^o(t, x)]} &&&& X \\
	\\
	&& {\Omega'} &&&& {\Omega'} \\
	\\
	Y &&&& {\Sum_{y:Y} p[\psi^o(t, y)]} &&&& Y
	\arrow["{\theta'(t)}"', from=3-3, to=3-7]
	\arrow["{\vartheta^o(t)^\ast \sigma}", from=1-1, to=1-5]
	\arrow["{\vartheta^u(t)}", from=1-5, to=1-9]
	\arrow["{\psi^o(t)^\ast \sigma}"', from=5-1, to=5-5]
	\arrow["{\psi^u(t)}"', from=5-5, to=5-9]
	\arrow["f"', from=1-1, to=5-1]
	\arrow["f", from=1-9, to=5-9]
	\arrow["{\phi\circ\pi_Y}"', from=5-9, to=3-7]
	\arrow["{\phi\circ\pi_X}", from=1-9, to=3-7]
	\arrow["{\phi\circ\pi_Y}", from=5-1, to=3-3]
	\arrow["{\phi\circ\pi_X}"', from=1-1, to=3-3]
  \end{tikzcd}\]
  That is, \(\RDynT{T}(p, \phi)(\vartheta) := (\phi\circ\pi_X, \vartheta^o,
  \vartheta^u)\) and \(\RDynT{T}(p, \phi)(f)\) is given by the same underlying
  map \(f : X \to Y\) on state spaces.

  \begin{proof}
    We need to check: the naturality condition of Definition \ref{def:poly-rdyn}
    for both \(\RDynT{T}(\varphi, \theta)(\vartheta)\) and \(\RDynT{T}(p,
    \phi)(\vartheta)\); functoriality of \(\RDynT{T}(\varphi, \theta)\) and
    \(\RDynT{T}(p, \phi)\); and (pseudo)functoriality of \(\RDynT{T}\) with
    respect to both morphisms of polynomials and of closed measure-preserving
    systems.

    We begin by checking that the conditions of Definition \ref{def:poly-rdyn}
    are satisfied by the objects
    \[ \RDynT{T}(\varphi, \theta)(\pi_X, \vartheta^o, \vartheta^u) : \RDynT{T}(q,
    \theta)\] and morphisms \[\RDynT{T}(\varphi, \theta)(f) : \RDynT{T}(\varphi,
    \theta)(\pi_X, \vartheta^o, \vartheta^u) \to \RDynT{T}(\varphi,
    \theta)(\pi_Y, \psi^o, \psi^u)\] in the image of \(\RDynT{T}(\varphi,
    \theta)\). Given a section \(\tau : q(1) \to
    \Sum_{j:q(1)} q[j]\) of q, we need to check that the closure
    \(\RDynT{T}(\varphi, \theta)(\vartheta)^\tau\) forms a closed random
    dynamical system in \(\Cat{Cat}(\deloop{\Tt}, \cat{E})/\theta\). That is to
    say, for all \(t : \Tt\) and sections \(\tau\), we need to check that the
    following naturality square commutes:
    
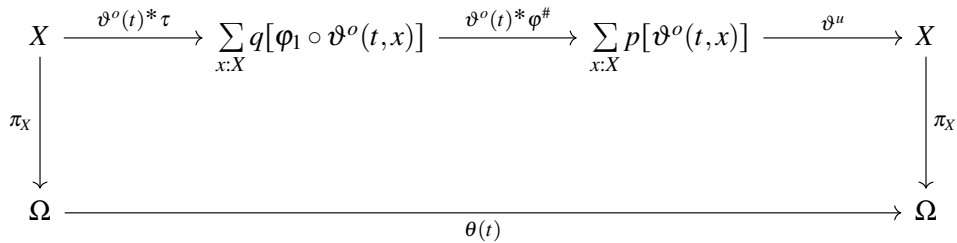
\[\begin{tikzcd}
	X && {\Sum_{x:X} q[\varphi_1 \circ \vartheta^o(t,x)]} && {\Sum_{x:X} p[\vartheta^o(t,x)]} && X \\
	\\
	\Omega &&&&&& \Omega
	\arrow["{\vartheta^o(t)^\ast \tau}", from=1-1, to=1-3]
	\arrow["{\vartheta^o(t)^\ast \varphi^\#}", from=1-3, to=1-5]
	\arrow["{\vartheta^u}", from=1-5, to=1-7]
	\arrow["{\pi_X}"', from=1-1, to=3-1]
	\arrow["{\theta(t)}"', from=3-1, to=3-7]
	\arrow["{\pi_X}", from=1-7, to=3-7]
    \end{tikzcd}\]
    As before, we find that \(\varphi^\# \circ \varphi_1^\ast \tau\) is a
    section of \(p\), so that commutativity of the diagram above reduces to
    commutativity of the diagram in Definition \ref{def:poly-rdyn}. Similarly,
    given a morphism \(f : (\pi_X, \vartheta^o, \vartheta^u) \to (\pi_Y, \psi^o,
    \psi^u)\), we need to check that the diagram in Proposition
    \ref{prop:poly-rdyn} induced for \(\RDynT{T}(\varphi, \theta)(f)\) commutes
    for all times \(t : \Tt\) and sections \(\tau\) of \(q\). But given such a section \(\tau\),
    the diagram for \(\RDynT{T}(\varphi, \theta)(f)\) reduces to that for \(f\)
    and the section \(\varphi^\# \circ \varphi_1^\ast \tau\) of \(p\), which
    commutes \textit{ex hypothesi}; and functoriality of \(\RDynT{T}(\varphi,
    \theta)\) follows immediately.

    Next, we check that the conditions of Definition \ref{def:poly-rdyn} are
    satisfied in the image of \(\RDynT{T}(p, \phi)\). It is clear by the
    definition of the action of \(\RDynT{T}(p, \phi)\) that the condition that
    the diagram in Proposition \ref{prop:poly-dyn-idx} commutes is satisfied,
    from which it follows by pasting that \(\RDynT{T}(p, \phi)\) is
    functorial. We therefore just have to check the induced diagram in
    Definition \ref{def:poly-rdyn} commutes. Consider the following diagram:
    \[\begin{tikzcd}
	X && {\Sum_{x:X} p[\vartheta^o(t, x)]} && X \\
	\\
	\Omega &&&& \Omega \\
	\\
	{\Omega'} &&&& {\Omega'}
	\arrow["{\pi_X}"', from=1-1, to=3-1]
	\arrow["{\vartheta^o(t)^\ast \sigma}", from=1-1, to=1-3]
	\arrow["{\vartheta^u(t)}", from=1-3, to=1-5]
	\arrow["{\pi_X}", from=1-5, to=3-5]
	\arrow["{\theta(t)}", from=3-1, to=3-5]
	\arrow["\phi"', from=3-1, to=5-1]
	\arrow["\phi", from=3-5, to=5-5]
	\arrow["{\theta'(t)}", from=5-1, to=5-5]
    \end{tikzcd}\]
    The top square commutes \textit{ex hypothesi}, the bottom square commutes by
    the definition of morphism of closed measure-preserving dynamical systems
    (Proposition \ref{prop:metric-sys}), and the outer square is the induced
    diagram we need to check, which therefore commutes by the pasting of
    commuting squares.

    Finally, we check that \(\RDynT{T}\) is functorial with respect to morphisms
    of polynomials and morphisms of closed measure-preserving dynamical systems.
    These reduce to
    checking that pullback and composition are functorial, which we again leave
    to the dedicated reader.
  \end{proof}
\end{prop}

\begin{defn}[Open bundle dynamical system] \label{defn:bdyn-pbth}
  Let \(p, b : \Poly{E}\) be polynomials in \(\cat{E}\), and let \(\theta :=
  (\theta(\ast), \theta^o, \theta^u) : \DynT{T}(b)\) be an open dynamical system
  over \(b\).  An \textnormal{open bundle dynamical system} over \((p, b, \theta)\)
  is a pair \((\pi_{\vartheta\theta}, \vartheta)\) where \(\vartheta :=
  (\vartheta(\ast), \vartheta^o, \vartheta^u) : \DynT{T}(p)\) is an open
  dynamical system over \(p\) and \(\pi_{\vartheta\theta} : \vartheta(\ast) \to
  \theta(\ast)\) is a bundle in \(\cat{E}\), such that, for all time \(t : \Tt\)
  and sections \(\sigma\) of \(p\) and \(\varsigma\) of \(b\), the following
  diagrams commute, thereby inducing a bundle of closed dynamical systems
  \(\pi^{\sigma\varsigma}_{\vartheta\theta} : \vartheta^\sigma \to
  \theta^\varsigma\) in \(\Cat{Cat}(\deloop{\Tt}, \cat{E})\):
  \[\begin{tikzcd}
	\vartheta(\ast) && {\Sum_{w:\vartheta(\ast)} p[\vartheta^o(t, w)]} && {\vartheta(\ast)} \\
	\\
	\theta(\ast) && {\Sum_{x:\theta(\ast)} b[\theta^o(t,x)]} && \theta(\ast)
	\arrow["{\pi_{\vartheta\theta}}", from=1-1, to=3-1]
	\arrow["{\pi_{\vartheta\theta}}", from=1-5, to=3-5]
	\arrow["{\vartheta^o(t)^\ast \sigma}", from=1-1, to=1-3]
	\arrow["{\vartheta^u(t)}", from=1-3, to=1-5]
	\arrow["{\theta^o(t)^\ast \varsigma}", from=3-1, to=3-3]
	\arrow["{\theta^u(t)}", from=3-3, to=3-5]
  \end{tikzcd}\]
\end{defn}

\begin{prop}[Category of open bundle dynamical systems over \((p,b)\)] \label{prop:bdyn-pbth-cat}
  Let \(p, b : \Poly{E}\) be polynomials in \(\cat{E}\), and let \(\theta :=
  (\theta(\ast), \theta^o, \theta^u) : \DynT{T}(b)\) be an open dynamical system
  over \(b\).  Open bundle dynamical systems over \((p, b, \theta)\) form the
  objects of a category \(\BDynT{T}(p, b, \theta)\). Morphisms \(f :
  (\pi_{\vartheta\theta}, \vartheta) \to (\pi_{\varrho\theta}, \varrho)\) are
  maps \(f : \vartheta(\ast) \to \varrho(\ast)\) in \(\cat{E}\) making the
  following diagram commute for all times \(t : \Tt\) and sections \(\sigma\) of
  \(p\) and \(\varsigma\) of \(b\):
  \[\begin{tikzcd}
	{\vartheta(\ast)} &&&& {\Sum_{w:\vartheta(\ast)} p[\vartheta^o(t, w)]} &&&& {\vartheta(\ast)} \\
	\\
	&& {\theta(\ast)} && {\Sum_{x:\theta(\ast)} b[\theta^o(t,x)]} && {\theta(\ast)} \\
	\\
	{\varrho(\ast)} &&&& {\Sum_{y:\varrho(\ast)} p[\varrho^o(t, y)]} &&&& {\varrho(\ast)}
	\arrow["{\vartheta^o(t)^\ast \sigma}", from=1-1, to=1-5]
	\arrow["{\vartheta^u(t)}", from=1-5, to=1-9]
	\arrow["{\theta^o(t)^\ast \varsigma}", from=3-3, to=3-5]
	\arrow["{\theta^u(t)}", from=3-5, to=3-7]
	\arrow["{\pi_{\vartheta\theta}}"', from=1-1, to=3-3]
	\arrow["{\pi_{\varrho\theta}}", from=5-1, to=3-3]
	\arrow["{\pi_{\vartheta\theta}}", from=1-9, to=3-7]
	\arrow["{\pi_{\varrho\theta}}"', from=5-9, to=3-7]
	\arrow["{\varrho^o(t)^\ast \sigma}"', from=5-1, to=5-5]
	\arrow["{\varrho^u(t)}"', from=5-5, to=5-9]
	\arrow["{f}"', from=1-1, to=5-1]
	\arrow["{f}", from=1-9, to=5-9]
  \end{tikzcd}\]
  That is, \(f\) is a map on the state spaces that induces a morphism
  \((\pi_{\vartheta\theta}, \vartheta^\sigma) \to (\pi_{\varrho\theta},
  \varrho^\sigma)\) in \(\Cat{Cat}(\deloop{\Tt}, \cat{E})/\theta^\varsigma\) of
  bundles of the closures. Identity morphisms are the corresponding identity
  maps, and composition is by pasting.
\end{prop}

\begin{prop}[Opindexed category of open bundle dynamical systems] \label{prop:bdyn-bth-idx}
  Varying the polynomials \(p\) in \(\BDynT{T}(p, b, \theta)\) induces an
  opindexed category \(\BDynT{T}({-}, b, \theta) : \Poly{E} \to \Cat{Cat}\). On
  polynomials \(p\), it returns the categories \(\BDynT{T}(p, b, \theta)\) of
  Proposition \ref{prop:bdyn-pbth-cat}. On morphisms \(\varphi : p \to q\) of
  polynomials, define the functors \(\BDynT{T}(\varphi, b, \theta) :
  \BDynT{T}(p, b, \theta) \to \BDynT{T}(q, b, \theta)\) as in Proposition
  \ref{prop:poly-rdyn-idx}. That is, suppose \((\pi_{\vartheta\theta},
  \vartheta) : \BDynT{T}(p, b, \theta)\) is object (open bundle dynamical
  system) in \(\BDynT{T}(p, b, \theta)\), where \(\vartheta := (\vartheta(\ast),
  \vartheta^o, \vartheta^u)\). Then its image \(\BDynT{T}(\varphi, b,
  \theta)(\pi_{\vartheta\theta}, \vartheta)\) is defined as the pair
  \((\pi_{\vartheta\theta}, \varphi\vartheta)\), where \(\varphi\vartheta :=
  (\vartheta(\ast), \phi_1 \circ \vartheta^o, \vartheta^u \circ
    {\vartheta^o}^\ast \varphi^\#)\). On morphisms \(f : (\pi_{\vartheta\theta},
    \vartheta) \to (\pi_{\varrho\theta}, \varrho)\), \(\BDynT{T}(\varphi, b,
    \theta)(f)\) is again given by the same underlying map \(f : \vartheta(\ast)
    \to \varrho(\ast)\) of state spaces.
  \begin{proof}
    The proof amounts to the proof for Proposition \ref{prop:poly-rdyn-idx} that
    \(\RDynT{T}(\varphi, \theta)\) constitutes an indexed category, except that
    the closed base dynamical system \(\theta\) of that Proposition is here
    replaced, for any section \(\varsigma\) of \(b\), by the closure
    \(\theta^\varsigma\) by \(\varsigma\) of the open dynamical system \(\theta
    : \DynT{T}(b)\) of the present Proposition. The proof goes through
    accordingly, since the relevant diagrams are guaranteed to commute for any
    such \(\varsigma\) by the conditions in Definition \ref{defn:bdyn-pbth} and
    Proposition \ref{prop:bdyn-pbth-cat}.
  \end{proof}
\end{prop}

\begin{prop}[Doubly-opindexed category of open bundle dynamical systems] \label{prop:bdyn-b-idx}
  Letting the base system \(\theta\) also vary induces a doubly-opindexed category
  \(\BDynT{T}({-}, b, {=}) : \Poly{E} \times \DynT{T}(b) \to \Cat{Cat}\). Given
  a polynomial \(p : \Poly{E}\) and morphism \(\phi : \theta \to \rho\) in
  \(\DynT{T}(b)\), the functor \(\BDynT{T}(p, b, \phi) : \BDynT{T}(p, b, \theta)
  \to \BDynT{T}(p, b, \rho)\) is defined by post-composition, as in Proposition
  \ref{prop:poly-rdyn-idx} for the action of \(\RDynT{T}\) on morphisms of the
  base systems there. More explicitly, such a morphism \(\phi\) corresponds to a
  map \(\phi : \theta(\ast) \to \rho(\ast)\) of state spaces in
  \(\cat{E}\). Given an object \((\pi_{\vartheta\theta}, \vartheta)\) of
  \(\BDynT{T}(p, b, \theta)\), we define \(\BDynT{T}(p, b,
  \phi)(\pi_{\vartheta\theta}, \vartheta) := (\phi \circ \pi_{\vartheta\theta},
  \vartheta)\). Given a morphism \(f : (\pi_{\vartheta\theta}, \vartheta) \to
  (\pi_{\varrho\theta}, \varrho)\) in \(\BDynT{T}(p, b, \theta)\), its image
  \(\BDynT{T}(p, b, \phi)(f) : (\phi\circ\pi_{\vartheta\theta}, \vartheta) \to
  (\phi\circ\pi_{\varrho\theta}, \varrho)\) is given by the same underlying map
  \(f : \vartheta(\ast) \to \varrho(\ast)\) of state spaces.
  \begin{proof}
    As for Proposition \ref{prop:bdyn-bth-idx}, the proof here amounts to the
    proof for Proposition \ref{prop:poly-rdyn-idx} that \(\RDynT{T}(p, \phi)\)
    constitutes an indexed category, except again the closed systems are
    replaced by (the appropriate closures of) open ones, and the
    measure-preserving structure is forgotten.
  \end{proof}
\end{prop}

\begin{nproof}[Proof of the Laplace approximation] \label{proof:laplace-approx}
First note that the KL divergence is bounded from above by the free energy since $\log p_{\gamma\klcirc\pi}(y)$ is always negative.

Next, we can write the density functions as:
\begin{gather*}
\log p_\gamma (y | x) = \frac{1}{2} \innerprod{\epsilon_\gamma}{{\Sigma_\gamma}^{-1} {\epsilon_\gamma}} - \log \sqrt{(2 \pi)^{|Y|} \det \Sigma_\gamma } \\
\log p_{\rho_\pi} (x | y) = \frac{1}{2} \innerprod{\epsilon_{\rho_\pi}}{{\Sigma_{\rho_\pi}}^{-1} {\epsilon_{\rho_\pi}}} - \log \sqrt{(2 \pi)^{|X|} \det \Sigma_{\rho_\pi} } \\
\log p_\pi (x) = \frac{1}{2} \innerprod{\epsilon_\pi}{{\Sigma_\pi}^{-1} {\epsilon_\pi}} - \log \sqrt{(2 \pi)^{|X|} \det \Sigma_\pi }
\end{gather*}
where for clarity we have omitted the dependence of \(\Sigma_\gamma\) on \(x\) and \(\Sigma_{\rho_\pi}\) on \(y\), and where
\begin{align*}
\epsilon_\gamma : Y \times X \to Y & : (y, x) \mapsto y - \mu_\gamma(x) \, , \\
\epsilon_{\rho_\pi} : X \times Y \to X & : (x, y) \mapsto x - \mu_{\rho_\pi}(y) \, , \\
\epsilon_\pi : X \times 1 \to X & : (x, \ast) \mapsto x - \mu_\pi \, .
\end{align*}
Then, note that we can write the free energy \(\Fa(y)\) as the difference between expected energy and entropy:
\begin{align*}
  \Fa(y)
  &= \E_{x \sim \rho_\pi(y)} \left[ \log \frac{p_{\rho_\pi}(x|y)}{p_\gamma(y|x) \cdot p_\pi(x)} \right] \\
  &= \E_{x \sim \rho_\pi(y)} \left[ - \log p_\gamma(y|x) - \log p_\pi (x) \right]
     - S_X \left[ \rho_\pi(y) \right] \\
  &= \E_{x \sim \rho_\pi(y)} \left[ E_{(\pi,\gamma)}(x,y) \right] - S_X \left[ \rho_\pi(y) \right]
\end{align*}
Next, since the eigenvalues of \(\Sigma_{\rho_\pi}(y)\) are small for all \(y : Y\), we can approximate the expected energy by its second-order Taylor expansion around the mean \(\mu_{\rho_\pi}(y)\):
\begin{align*}
\Fa(y) \approx & \;
E_{(\pi,\gamma)}(\mu_{\rho_\pi}(y), y)
+ \frac{1}{2} \innerprod{\epsilon_{\rho_\pi}\left( \mu_{\rho_\pi}(y), y\right)}{\left(\partial_x^2 E_{(\pi,\gamma)}\right)\left( \mu_{\rho_\pi}(y), y\right) \cdot \epsilon_{\rho_\pi}\left( \mu_{\rho_\pi}(y), y\right)} \\
& \;\, - S_X \big[ \rho_\pi(y) \big] \, .
\end{align*}
where \(\left(\partial_x^2 E_{(\pi,\gamma)}\right)\left( \mu_{\rho_\pi}(y), y\right)\) is the Hessian of \(E_{(\pi,\gamma)}\) with respect to \(x\) evaluated at \((\mu_{\rho_\pi}(y), y)\).

Note that
\begin{equation} \label{eq:laplace-trace-sigma}
\innerprod{\epsilon_{\rho_\pi}\left( \mu_{\rho_\pi}(y), y\right)}{\left(\partial_x^2 E_{(\pi,\gamma)}\right)\left( \mu_{\rho_\pi}(y), y\right) \cdot \epsilon_{\rho_\pi}\left( \mu_{\rho_\pi}(y), y\right)}
=
\tr \left[ \left(\partial_x^2 E_{(\pi,\gamma)}\right)\left( \mu_{\rho_\pi}(y), y\right) \, \Sigma_{\rho_\pi}(y) \right] \, ,
\end{equation}
that the entropy of a Gaussian measure depends only on its covariance,
\[
S_X \big[ \rho_\pi(y) \big]
= \frac{1}{2} \log \det \left( 2 \pi \, e \, \Sigma_{\rho_\pi}(y) \right) \, ,
\]
and that the energy \(E_{(\pi,\gamma)}(\mu_{\rho_\pi}(y), y)\) does not depend on \(\Sigma_{\rho_\pi}(y)\). We can therefore write down directly the covariance \(\Sigma_{\rho_\pi}^\ast(y)\) minimizing \(\Fa(y)\) as a function of \(y\). We have
\[
\partial_{\Sigma_{\rho_\pi}} \Fa(y) \approx
\frac{1}{2} \left(\partial_x^2 E_{(\pi,\gamma)}\right)\left( \mu_{\rho_\pi}(y), y\right)
+ \frac{1}{2} {\Sigma_{\rho_\pi}}^{-1} \, .
\]
Setting \(\partial_{\Sigma_{\rho_\pi}} \Fa(y) = 0\), we find the optimum
\[
\Sigma_{\rho_\pi}^\ast (y) = \left(\partial_x^2 E_{(\pi,\gamma)}\right)\left( \mu_{\rho_\pi}(y), y\right)^{-1} \, .
\]
Finally, on substituting \(\Sigma_{\rho_\pi}^\ast (y)\) in equation \eqref{eq:laplace-trace-sigma}, we obtain the desired expression
\[
\Fa(y) \approx E_{(\pi,\gamma)}\left(\mu_{\rho_\pi}(y), y\right) - S_X \left[ \rho_\pi(y) \right] =: \Fa^L(y) \, .
\]
\qed\end{nproof}

\end{document}